\documentclass[12pt]{article}

\usepackage[utf8]{inputenc}
\usepackage[T1]{fontenc}
\usepackage{bbm}
\usepackage{stmaryrd}
\usepackage{latexsym}
\usepackage{amsfonts}
\usepackage{amsmath,amssymb,amscd,textcomp}
\usepackage{yfonts}
\usepackage{dsfont}
\usepackage{mathabx}
\usepackage[dvips]{graphicx}
\usepackage{epsfig}
\usepackage{psfrag}
\usepackage[font={small}]{caption}
\usepackage{color}
\usepackage{float}
\usepackage{upgreek}
\usepackage{tikz}
\usepackage{tikz-cd}
\usepackage{relsize}
\usepackage{bbold}
\usepackage{hyperref}
\usepackage{fullpage}
\usepackage{array}

\DeclareFontFamily{U}{txsyc}{}
\DeclareFontShape{U}{txsyc}{m}{n}{
   <-> txsyc%
}{}
\DeclareFontShape{U}{txsyc}{bx}{n}{
   <-> txbsyc%
}{}
\DeclareFontShape{U}{txsyc}{l}{n}{<->ssub * txsyc/m/n}{}
\DeclareFontShape{U}{txsyc}{b}{n}{<->ssub * txsyc/bx/n}{}
\DeclareSymbolFont{symbolsC}{U}{txsyc}{m}{n}
\SetSymbolFont{symbolsC}{bold}{U}{txsyc}{bx}{n}
\DeclareFontSubstitution{U}{txsyc}{m}{n}
\DeclareMathSymbol{\df}{\mathrel}{symbolsC}{"42}
\DeclareMathSymbol{\fd}{\mathrel}{symbolsC}{"43}
\DeclareMathSymbol{\lJoin}{\mathrel}{symbolsC}{"58}
\DeclareMathSymbol{\rJoin}{\mathrel}{symbolsC}{"59}

\newcommand{\cB}{\mathcal{B}}

\newcommand{\cL}{\mathcal{L}}

\newcommand{\cU}{\mathcal{U}}

\newcommand{\EE}{\mathbb{E}}

\newcommand{\PP}{\mathbb{P}}

\newcommand{\fs}{\mathfrak{s}}

\newcommand{\lt}{\left}
\newcommand{\me}{\medskip}

\newcommand{\rt}{\right}
\newcommand{\sm}{\smallskip}

\newcommand{\wi}{\widetilde}

\DeclareMathOperator*{\esssup}{ess\,sup}

\newcommand{\fo}{\forall\ }

\newcommand{\ind}[1]{\mathbbm{\un}_{\! #1}}
\newcommand{\un}{\mathds{1}}
\newcommand{\Var}{\mathrm{Var}}

\newcommand{\bq}{\begin{eqnarray*}}
\newcommand{\bqn}[1]{\begin{eqnarray}\label{#1}}
\newcommand{\eq}{\end{eqnarray*}}
\newcommand{\eqn}{\end{eqnarray}}
\newcommand{\wwtbp}{\par\hfill $\blacksquare$\par\me\noindent}

\newcommand{\thistitlepagestyle}{}

\newcommand{\ttsim}{\raise.17ex\hbox{$\scriptstyle\mathtt{\sim}$}}

\newcommand{\kh}{\kern .08em}

\newtheorem{pro}{Proposition} 
\newtheorem{cor}[pro]{Corollary}
\newtheorem{lem}[pro]{Lemma}
\newtheorem{theo}[pro]{Theorem}
\newtheorem{defi}[pro]{Definition}
\newtheorem{Rem}[pro]{Remark}

\renewcommand{\thepro}{\arabic{pro}}

\newcommand{\proof}{\par\me\noindent\textbf{Proof}  \par\sm\noindent }
\newcommand{\prooff}[1]{\par\me\noindent\textbf{#1}\par\sm\noindent}

\newcommand{\comment}[1]{}

\allowdisplaybreaks

\newcommand{\Hm}[1]{\leavevmode{\marginpar{\tiny%
$\hbox to 0mm{\hspace*{-0.5mm}$\leftarrow$\hss}%
\vcenter{\vrule depth 0.1mm height 0.1mm width \the\marginparwidth}%
\hbox to 0mm{\hss$\rightarrow$\hspace*{-0.5mm}}$\\\relax\raggedright
#1}}}
\usepackage[notcite,notref]{
}

\usepackage{xcolor}

\newcommand{\eps}{\epsilon}
\title{Cutoff in separation for compact harmonic manifold}
 \author{Magalie Bénéfice${}^{(1)}$, Kol\'eh\`e Coulibaly-Pasquier${}^{(1)}$ }

 \date{\vbox{\copy0
 \copy1
 \copy2
}
 }

\begin{document}


\setbox1=\vbox{
\large
\begin{center}
${}^{(1)}$ Institut \'Elie Cartan de Lorraine, UMR 7502\\
Universit\'e de Lorraine and CNRS
\end{center}
} 

\setbox4=\vbox{
\hbox{magalie.benefice@univ-lorraine.fr\\[1mm]}
\hbox{Institut \'Elie Cartan de Lorraine\\}
\hbox{Universit\'e de Lorraine}
}
\setbox5=\vbox{
\hbox{kolehe.coulibaly@univ-lorraine.fr\\[1mm]}
\hbox{Institut \'Elie Cartan de Lorraine\\}
\hbox{Universit\'e de Lorraine}
}


\maketitle
\thistitlepagestyle
\abstract{We show that cutoff in separation occurs for Brownian motion in some families of compact harmonic manifolds. We compute the cutoff time and windows in four families of compact harmonic manifold namely $\mathbb{S}^n $, $\mathbb{CP}^n$, $\mathbb{HP}^n$ and  $\mathbb{RP}^n$ (the first three families are the only  families of compact simply connected harmonic manifolds, see \cite{Szabo}).
The proof is based on sharp strong stationary times and sufficiently accurate asymptotic expansions of their means and variances.}

\vfill\null
{\small
\par\vskip.3cm
}\par


\section{Overview}
\subsection{Introduction}
Consider a Markov process $X$. Suppose that it admits a unique invariant distribution $\mu$ and that $X$ converges in distribution to $\mu$ in large times. It is then natural to try to quantify this convergence. This is particularly relevant when the cardinal or the dimension $d$ of the state space of $X$ is considered to be large. In these cases a cutoff phenomenon, i.e., an abrupt transition to equilibrium, can sometimes be observed and the time for this transition as well as its abruptness (window) can be studied. We are here interested in these phenomenons for $X$ being a Brownian motion on the Riemannian manifold $M$, $M$ being the high dimensional sphere or a high dimensional projective space $\mathbb{KP}^n$ with $\mathbb{K}\in\{\mathbb{R},\mathbb{C},\mathbb{H}\}$. These manifolds are compact symmetric spaces of rank one and they are  harmonic manifolds, i.e., their volume forms and equivalently the mean curvature of spheres are radial around any given point $x_0$.
 The cutoff properties of Brownian motions on these spaces, relative to the total variation distance and the $L^p$-norms ($p\in[1,\infty]$) have been well studied and extended to other compact symmetric spaces by \cite{MR1306030, SCC, MR3201989}.
 More recently, \cite{Salez} produced a criteria for  cutoff in total variation, for non negatively curved diffusions, when an additional product condition hypothesis is satisfied, using entropy and varentropy. In the present paper we consider cutoff relative to another  way to measure the gap between the distribution of  $X_t$ and its invariant measure: the separation distance. Even if total variation distance, separation distance and $L^{\infty}$-norm can be compared, there is no equivalence for the associated cutoff in general~\cite{zbMATH06618510}. The method used in this article is then quite different than for \cite{MR1306030,SCC,MR3201989, Salez}. We are using two tools: strong stationary times and dual processes.

Strong stationary times where first introduced in~\cite{AD} by Aldous and Diaconis in the context of card shuffling.  A strong stationary time is a stopping time $\tau$ for $X$, independent of $X_{\tau}$ and such that $X_{\tau} \overunderset{\mathcal{L}}{}\sim \mu$. In particular it provides, in the discrete case, upper-bounds for separation distance $\mathfrak{s}$
(see~\cite{ADbis}):
\begin{equation}\label{EsT} 
\mathfrak{s}(\mathcal{L}(X_t),\mu)\leq \PP(\tau > t) ,\, \forall t>0.
\end{equation}
In~\cite{ADbis}, Aldous and Diaconis proved that there always exists at last one strong stationary time such that \eqref{EsT} becomes an equality. Such stopping times are called sharp strong stationary times. Even if  sharp  strong stationary times are well-adapted to the study of separation cutoffs, they are not evident to construct and to study  in general.

For the fundamental case of top at random card shuffling~(described in~\cite{AD}), where $X$ is a Markov chain in the permutations, their construction is elegant and quite straightforward. The constructed sharp strong stationary time can in fact be seen as the absorption time of a dual Markov process $D$ that takes values in the power set of the permutations. Here the absorbing state of $D$ is the full state space of $X$. In the same spirit, Diaconis and Fill~\cite{MR1071805} proposed the study of strong stationary times in the case of various Markov chains via the construction of dual processes intertwined with $X$ and the study of their absorption times. 
The theory has been extended by Miclo~\cite{M19} to study diffusion processes in $\mathbb{R}$ by considering dual diffusion processes whose values are intervalls.
In \cite{zbMATH07470497,arnaudon:hal-03037469} the authors considered the generalization of the theory to a Riemannian manifold $M$. In \cite{zbMATH07470497} the authors considered the construction and depiction of a candidate for a dual diffusion process $D$ whose values are domains of $M$. In~\cite{arnaudon:hal-03037469} the authors constructed several couplings  $(X_t,D_t)_{t \ge 0}$  to satisfy the intertwined relations. 
To get exactly an intertwined relation in the sense of Diaconis and Fill and obtain strong or sharp strong stationary times, one remaining condition is to be able to start $D$ from the singleton $\{X_0\}$. The other one is that the absorption state of $D$ must be $M$, and the absorption time have to be a.s. finite. 

In particular, the first condition is verified when $M$ satisfies some radial properties around $X_0$. The dual process $D$ can then be chosen taking values into geodesic balls centred in $X_0$ and can be described by its radius, bringing back the problem to the study of a first hitting time for the radius $(R_t)_{t \ge 0}$. The case of sphere and rotationally symmetric manifolds (where the radial properties are satisfied at last around one point) have been studied in~\cite{sphere,ACM}. Note that projective spaces are not rotationally symmetric manifolds in the sense of~\cite{ACM}.

In this paper we use the radial  properties of harmonic manifolds to extend the method to high dimensional projective spaces and improve the results for the sphere (especially the window). By studying both expectation and variance for the absorption time, we obtain cutoff in separation results for Brownian motion with  windows. For the variance, we use a method inspired from the stationary phase approximation up to a parametrization to obtain equivalents of multiple integrals.

Note that, in particular, the global method does not directly use the knowledge of the spectrum of the Laplace Beltrami operator, contrary to the  precedent result  concerning cutoff for  $L^p $ distance e.g. \cite{MR1306030, SCC, MR3201989} .

\subsection{Main Results}

Set a family $(M_n)_{n \in \mathbb{N}}$ of compact Riemannian manifolds with a dimension $(d_n)_{n \in \mathbb{N}}$. For all $n$, denote by $(X^n_t)_{t\geq 0}$ the Brownian motion on $M_n$ and by $\mathcal{U}_n$ its invariant distribution, namely the uniform distribution on $M_n$. By ergodicity, for $t\to\infty$, $X^n_t$ converges in distribution to $\mathcal{U}_n$. Our aim is to study the behaviour of this convergence for $n$ large. For this, we are using the distance in separation to evaluate the discrepancy between the distributions $\mathcal{L}(X^n_t)$ and $\mathcal{U}_n$:
\begin{defi}
For $\mu$ and $\nu$ two probability measures defined on the same measurable space, the separation distance $\mathfrak{s}$ between $\mu$ and $\nu$ is given by:
$$\mathfrak{s}(\mu,\nu):=\underset{\nu}{\esssup}\left(1-\frac{d\mu}{d\nu}\right)$$ where 
$d\mu/d\nu$ is the Radon-Nikodym derivative of the absolutely continuous part  of $\mu$ with respect to $\nu$ and
$\underset{\nu}{\esssup}$ is the essential supremum with respect to $\nu$.

\end{defi}
\begin{Rem}
The separation distance takes values in $[0,1]$.
Note that  $\mathfrak{s}$ is not a distance since it is not symmetric in its arguments.
\end{Rem}

 Various definitions for the cutoff phenomenons and their windows can be found in Definition 2.1. from~\cite{SCC}. In a similar way, in the case of the separation distance we will use the following definition:
 \begin{defi}
 The family $(X^n)_n$ of Brownian motions on the growing dimensional family $(M_n)_n$ admits a  cutoff in separation with cutoff time $(a_n)_{n}$, and window $(b_n)_{n}$ 
    if:     $b_n =o(a_n) $ and
 \begin{itemize}
    \item $\lim\limits_{r\to \infty}\limsup \limits_{n\to\infty}\mathfrak{s}(\mathcal{L}(X^n_{a_n+rb_n}),\mathcal{U}_n)=0$
    \item $\lim\limits_{r\to \infty}\liminf \limits_{n\to\infty}\mathfrak{s}(\mathcal{L}(X^n_{a_n-rb_n}),\mathcal{U}_n)= 1.$
\end{itemize}

 \end{defi} 
 Our main result is about cutoff in separation results for the high dimensional spheres, complex projective spaces and quaternionic projective spaces, i.e., compact simply connected harmonic manifolds. We also obtain these results for real projective spaces. 
 \begin{theo}\label{Tmain}
Denote by $M_n$ one of the following families of Riemannian manifolds (with the usual Riemannian metric): 
\begin{itemize}
    \item the sphere $\mathbb{S}^n$ of diameter $\pi$ and of $\mathbb R$-dimension $d_n=n$;
    \item the projective space $\mathbb{KP}^n$ of diameter $\pi$ with $\mathbb{K}\in\{\mathbb{R}, \mathbb{C}, \mathbb{H}\}$ and of $\mathbb R$-dimension $d_n=an$, where $a=dim_{\mathbb{R}}(\mathbb{K})$.
\end{itemize}  The (two-time accelerated) Brownian motion $(X^n)_n$ on $M_n$, admits a cutoff in separation with cutoff time $a_n\sim\begin{cases} \frac{\ln(d_n)}{d_n}\text{ if }M_n=\mathbb{S}^n\\
\frac{2\ln(d_n)}{d_n}\text{ if }M_n=\mathbb{KP}^n
\end{cases}$ and a  window  $b_n$ of order $\frac{1}{d_n}$.
More details about the cutoff  will be detailed in Theorems~\ref{T1}, \ref{T2}, \ref{T3}, \ref{T4} depending of the chosen families.
 \end{theo}
 \begin{Rem}
 In the case of the sphere, the above results improve the window $b_n=o\left(\frac{\sqrt{\ln(n)}}{n}\right)$ obtained in~\cite{sphere}.
 \end{Rem}
 \begin{Rem}
 If the diameter chosen for $M_n$ in Theorem \ref{Tmain} is $\alpha>0$, then for $(a_n)_n$ and $(b_n)_n$ as in Theorem \ref{Tmain}, there is a separation  cutoff phenomenon with cutoff time $\left(\frac{\alpha}{\pi}\right)^2 a_n$ and window $ b_n$.
 \end{Rem}
 \begin{Rem}
  Saloff-Coste and Chen~\cite{MR1306030,SCC} investigated these cases for cutoff with 
  \begin{itemize}
      \item total variation distance: $$d_{TV}(\mu,\nu):= \sup\limits_{A \text{ measurable}}|\mu(A)-\nu(A)|;$$
      \item for $L^p$-norms, for $1\leq p\leq \infty$: $$d_p(\mu,\nu)=\left\|\frac{d\mu}{d\nu}-1\right\|_{L^p(\nu)}.$$
      If $\mu$ is not absolutely continuous with respect to $\nu$, we take $d_1(\mu,\nu)=2$ and $d_p(\mu,\nu)=+\infty$ for $p> 1$.
  \end{itemize}
 Note that total variation cutoff and $L^1$ cutoff are equivalent as the two quantities only differ from a multiplicative constant. In general, even if $L^p$-norms and separation distance can be compared, there is no equivalence between the associated cutoff phenomenons in general (see~\cite{zbMATH06618510} for counter examples). Table~\ref{Table} sum up the obtained results, in this work we are interested in completing the last column (Theorem \ref{Tmain}):
  \begin{table}[h]
  \caption{cutoff for (two-time accelerated) Brownian motions on harmonic manifolds}
  \label{Table}
 
\begin{tabular}{l|m{4cm}|m{4cm}|m{4cm}|}
\cline{2-4}
                  &  \centering$d_{TV}$, $L^p$-norm & \centering$L^{\infty}$-norm & {\centering{Separation distance}}\\ \hline
\multicolumn{1}{|l|}{$\mathbb{S}^n$, $d_n=n$}   &$a_n\sim\frac{\ln(d_n)}{2d_n}$\newline $b_n$ of order $\frac{1}{d_n}$ optimal& $a_n\sim\frac{\ln(d_n)}{d_n}$\newline $b_n$ of order $\frac{1}{d_n}$ optimal  &  $a_n\sim\frac{\ln(d_n)}{d_n}$\newline $b_n$ of order $\frac{1}{d_n}$ \\[5pt]
\hline
\multicolumn{1}{|l|}{$\mathbb{CP}^n$, $d_n=2n$}   &$a_n\sim\frac{2\ln(d_n)}{d_n}$& $a_n\sim\frac{4\ln(d_n)}{d_n}$ & { $a_n\sim\frac{2\ln(d_n)}{d_n}$ }\newline$b_n$ of order $\frac{1}{d_n}$ \\[6pt]
\hline
\multicolumn{1}{|l|}{$\mathbb{HP}^n$, $d_n=4n$}   & $a_n\sim\frac{2\ln(d_n)}{d_n}$ & $a_n\sim\frac{4\ln(d_n)}{d_n}$ &  { $a_n\sim\frac{2\ln(d_n)}{d_n}$ }\newline$b_n$ of order $\frac{1}{d_n}$   \\[5pt]
\hline
\multicolumn{1}{|l|}{$\mathbb{RP}^n$, $d_n=n$}   & $a_n\sim\frac{2\ln(d_n)}{d_n}$ & $a_n\sim\frac{4\ln(d_n)}{d_n}$ & {$a_n\sim\frac{2\ln(d_n)}{d_n}$}\newline $b_n$ of order $\frac{1}{d_n}$  \\[5pt]
\hline
\end{tabular}
\end{table}
 \end{Rem}
 
\subsection{Paper structure}
In Section \ref{Sec0}, we remind the duality and intertwining theory studied in~\cite{zbMATH07470497} for Brownian motions on Riemannian manifolds and make a presentation of harmonic manifolds. In Section \ref{Sec0bis}, we apply the duality theory to spheres and projective spaces to construct sharp strong stationary times. In particular we give general formulae for the expectation and the variance of these sharp strong stationary times in function of the volume form. In the rest of the paper, for each considered manifold, we compute expansions or equivalences for the expectation and the variance of the sharp strong stationary time, then use it to obtain detailed cutoff results. Section \ref{Sec1} deals with the complexe projective space, Section \ref{Sec2} with the quaternionic projective space, Section \ref{Sec3} with the sphere and Section \ref{Sec4} with the real projective space.  
\section{Preliminaries}\label{Sec0}
\subsection{Dual process for Brownian motion}
Let $(M,g)$ be a $d $ dimensional Riemannian manifold,  $ \Delta $ the Laplace Beltrami operator, and $V_g$  the Riemannian volume measure.  Let $(X_t)_{t \ge 0}$ be the $ \Delta $ diffusion, namely the Brownian motion on $M$ accelerated by a factor $2$, with  invariant measure  $V_g$.  
Denote by $\mathcal{D}$ the set of non-empty compact and connected domains $D$ in $M$ such that its boundary $\partial D$ is smooth, as well as singletons of $M$.
 We are interested in studying diffusion processes with values in $\mathcal{D}$,
 with   infinitesimal generator $\mathfrak{L}$ (see~\cite{zbMATH07470497}), such that $\mathfrak{L}$ satisfies the following algebraic intertwining relation :
\begin{equation}\label{EIntertwing}
\mathfrak{L}\Lambda = \Lambda \Delta, 
\end{equation}
 where $\Lambda$ is the Markov kernel defined as follows: for all $ D \in \mathcal{D}$, i.e. domain or singleton  of $M$
\bq
  \fo A\in \cB(M),\qquad \Lambda(D,A)&\df &
\lt\{
\begin{array}{lcl}  
&\frac{ V_g (A\cap D)}{ V_g (D)} &\hbox{, if $\mu(D)>0$}\\
&\delta_x(A) & \hbox{, if $D=\{x\}$, with $x\in M$}.
\end{array}
\rt.
\eq

The general depiction of $\mathfrak{L}$ and of an associated diffusion process $(D_t)_{t \ge 0}$ has been done in~\cite{zbMATH07470497}, we say that $D$ is a  dual processes for $X$.

 Let us recall the  evolution equation for $D$ starting at a domain $D_0 $ :

\begin{equation}\label{Edomain}
dY_t=\left(\sqrt{2}dB_t+\left(2\frac{\bar{V}_g(\partial D_t)}{V_g(D_t)}-\sigma_{\partial D_t}(Y_t)\right)dt\right)N_{\partial D_t}(Y_t)\quad \forall Y_t\in\partial D_t
\end{equation}
where:
\begin{itemize}
\item  $\bar{V}_g$ is the  $(d-1)-$Hausdorff measure associated to $V_g $;
\item $(B_t)_{t \ge 0}$ is a one dimensional Brownian motion;
\item $N_{\partial D_t}(Y_t)$ is the exterior normal vector at point $Y_t$ on the boundary $\partial D_t$;
\item $\sigma_{\partial D_t}(Y_t)$ is the mean curvature at point $Y_t$ on the boundary $\partial D_t$.
\end{itemize}
Note that in general $D$ as a positive lifetime.

 In general the above equation \eqref{Edomain} is a partial stochastic differential equation, and is quite difficult to study, see for example~\cite{ACM1}.
 
\begin{Rem}\label{harm}
Nevertheless,  if in the manifold $M$, the mean curvature at any point on all geodesic spheres only depends on the radius, then the balls are preserved by \eqref{Edomain}, i.e. if Equation \eqref{Edomain} starts at a  {geodesic} ball $B(\wi 0,R(0))\df D_0$  then $ D_t $ remains a  {geodesic} ball.
It is then quite natural to apply this construction to manifolds satisfying this property, namely harmonic manifolds, see \eqref{ExprCourbure} below.
\end{Rem}

\subsection{Harmonic manifolds}
Consider $(M,g)$ a complete Riemannian manifold and $\tilde{0}\in M$. Around this point, the density of the Riemannian volume form $dV_g$ in polar coordinates can be described, at least for $r$ small enough, by:
$$\Theta_{\tilde{0}}(r,v)=\left|\det\left(T_{rv}\exp_{\tilde 0}\right)\right| {r^{n-1}}$$
where $v\in T_{\tilde 0}M$, $\|v\|=1$.

Harmonic manifolds are manifolds for which the volume form is radial. If such a property can be considered locally (see~\cite{Besse}), we will here restraint ourselves to examples where it is global:
\begin{defi}
$(M,g)$ is globally harmonic if there exists $\theta :\mathbb{R}^+\to\mathbb{R}$, such that $\Theta_{\tilde{0}}(r,v)=\theta(r)$ for all $r>0$, $v\in T_{\tilde 0}M$, $\|v\|=1$.
\end{defi}

\begin{Rem}
The name "harmonic" comes from the fact that, on globally harmonic manifolds, a function $f$ is harmonic (in the sense that $\Delta f=0$) if and only if the mean value on geodesic sphere is equal to its value at the center of the sphere, see 6.20, 6.21 in~\cite{Besse} .
\end{Rem}

Basic examples of harmonic manifolds are the two point homogeneous spaces, including the spheres as well as real, complexe and quaternionic projective spaces. The Lichnerowicz conjecture, as an attempt to classify harmonic manifolds, assumes that the converse statement is true. It is in fact false in general (see \cite{kreyssig2010introduction} for more details about the conjecture). Szab\'o proved in \cite{Szabo} that this is true for any compact simply connected harmonic manifolds. In this case, harmonic manifolds are  isometric either to the spheres $\mathbb{S}^n$, the projective spaces $\mathbb{KP}^n$ with $\mathbb{K}\in\{\mathbb{C}, \mathbb{H}\}$ or the octonionic projective space $\mathbb{0P}^2$.

As in Theorem \ref{Tmain}, we will consider $M_n$ being either the sphere $\mathbb{S}^n$ either one of the projective spaces $\mathbb{KP}^n$ of diameter $\pi$ with $\mathbb{K}\in\{\mathbb{R}, \mathbb{C}, \mathbb{H}\}$ studying as well the real projective space. 
\begin{pro}\label{PVHarmonic}
 For all $\tilde{0}\in M_n$, $v\in T_ {\tilde{0}}M_n$ such that $\|v\|=1$ and $r\in]0,\pi[$, the volume density at point $\exp_{\tilde 0}(rv)$ is given by:
\begin{equation}
\theta(r)=\begin{cases}
\sin^{n-1}(r) & \text{ if }M_n=\mathbb{S}^n\\
2^{an-1}\sin^{an-1}\left(\frac{r}{2}\right)\cos^{a-1}\left(\frac{r}{2}\right) & \text{ if }M_n=\mathbb{KP}^n,\, a=dim_{\mathbb{R}} (\mathbb{K})
\end{cases}.
\end{equation}
\end{pro}

\proof
For the sphere, the result is well known, obtained by using polar coordinates~\cite{Berger}.\\
The case of the projective spaces $\mathbb{KP}^n$ is usually considered with half of the present diameter and can be computed by using Jacobi fields.
 For $0\leq s\leq \pi$ denote by $\gamma(s):=\exp_{\tilde{0}}(sv)$ the geodesic defined by arc length starting from $\tilde 0$ in direction $v$. Set $(X_i)_{1\leq i\leq an}$ an orthonormal basis of $T_{\tilde 0}\mathbb{KP}^n$ such that $X_1=v$, and $\left(X_1(s),\hdots,X_{an}(s)\right)$ the orthonormal basis at $T_{\gamma(s)}\mathbb{KP}^n$ obtained by parallel transport along $\gamma$. From \cite{Berger}, $\Theta_{\tilde{0}}(r,v)$ can be obtained by considering the Jacobi fields $(Y_i)_{2\leq i\leq an}$ along $\gamma$ such that $Y_i(0)=0$ and $\dot Y_i(0)=X_i$ for all $2\leq i\leq an$ (here $\dot~$ denotes the covariant derivative of the vector field along $\gamma$):
 $$\Theta_{\tilde{0}}(r,v)=\det_{\{X_2(r),\hdots,X_{an}(r)\}}\left(Y_2(r),\hdots Y_{an}(r)\right).$$
 Such Jacobi fields can be constructed from 3.34 in~\cite{Besse}, up to a change of diameter and a renormalisation of some of the Jacobi fields: supposing that $(X_i)_{1\leq i\leq an}$ is an adapted basis (with respect to the projective structure, see 3.11 in \cite{Besse}), we have:
   $$Y_i(s)=\begin{cases} \sin\left(s\right)X_i(s) &\text{ if }2\leq i\leq a\\
   2\sin\left(\frac{s}{2}\right)X_i(s) &\text{ if }a+1\leq i\leq an\end{cases}$$
which gives the expected result.
\wwtbp

From Proposition \ref{PVHarmonic}, we can compute the mean curvature on $M_n$.

Indeed, for $0 < r < \pi$ and $v\in T_{\tilde 0}M$, $\|v\|=1$, the mean curvature at point $\exp_{\tilde 0}(rv)$ on the boundary of $B(\tilde 0,r)$ can be obtained from the formula (from~\cite{Szabo}):
\begin{equation}\label{ExprCourbure}
\sigma_{\partial B(\tilde 0,r)}\left(\exp_{\tilde 0}(rv)\right)=\frac{\partial_r\Theta_{\tilde 0}(r,v)}{\Theta_{\tilde 0}(r,v)}.\end{equation} 

In particular note that a complete Riemannian manifold $(M,g)$ is globally harmonic if and only if its mean curvature is constant on geodesic spheres.
\begin{pro}
 For all $\tilde{0}\in M_n$, $v\in T_ {\tilde{0}}M_n$ such that $\|v\|=1$ and $r\in]0,\pi[$, the mean curvature at point $\exp_{\tilde 0}(rv)$ on the boundary of $B(\tilde 0,r)$ is given by:
\begin{equation}\label{MCHarmonic}
\sigma(r)=\begin{cases}
(n-1)\cot(r) & \text{ if }M_n=\mathbb{S}^n\\
\frac{1}{2}\left((an-1)\cot\left(\frac{r}{2}\right)-(a-1)\tan\left(\frac{r}{2}\right)\right) & \text{ if }M_n=\mathbb{KP}^n,\, a=dim_{\mathbb{R}}(\mathbb{K})
\end{cases}.
\end{equation}
\end{pro}
\section{Strong stationary time on harmonic manifolds}\label{Sec0bis}
\subsection{Intertwined relations on harmonic manifolds}
We now return to consideration of dual processes of the Brownian motion. Using the radial properties of the volume form and of the mean curvature on harmonic manifolds as well as \eqref{ExprCourbure}, we get:
\begin{pro}
On harmonic manifolds $M_n $, let $\wi 0 \in M_n $, for any $R_0 \in ]0, \pi[$, and $D_0=B(\tilde 0,R_0)$, the process starting from $D_0$ and satisfying \eqref{Edomain} is $D_t=B(\tilde 0, R_t)$ where $(R_t)_{t \in [0, \tau_n]}$  satisfies 
$dR_t=\sqrt{2}dB_t+b(R_t)dt, 
$ and $\tau_n:=\inf\{t>0|\, R_t=\pi\}$,  with 
\begin{equation}\label{Ederive}
    b(r):=2\frac{\theta(r)}{\int_0^r\theta(r)dr}-\sigma(r)=2\frac{\theta(r)}{\int_0^r\theta(s)ds}-\frac{\theta'(r)}{\theta(r)}=\frac{d}{dr}\ln\left(\frac{\left(\int_0^r\theta(s)ds\right)^2}{\theta(r)}\right).
\end{equation}
\end{pro}

\proof
 Denote by $d_n$ the real dimension of the considered harmonic manifold $M_n $. Let $r\in]0,\pi[$,  recall that the cutlocus of $\wi 0 $ is equal to  $\text{cut}(\wi 0) = \{ x \in M_n,\text{ s.t. }\, d(\wi 0, x)= \pi \}$  e.g 3.35  in~\cite{Besse}, and as usual $V_g(\text{cut}(\wi 0)) =0 $. Hence the volume of a ball $V_g(B(\tilde 0, r)) = C_{d_n}\int_0^r\theta(s)ds$ where $ C_{d_n} $ is the volume of $ \mathbb{S}^{{d_n}-1}$, so $\bar{V}_g(\partial B(\tilde 0, r) ) = C_{d_n} \theta(r)$, and $ \sigma_{\partial B(\tilde 0,r)}= \frac{\theta'(r)}{\theta(r)} $.
Hence \eqref{Ederive} follow from the equation of \eqref{Edomain}.

\wwtbp
\begin{pro}\label{T01}
On $M_n$, the process starting from $\{\tilde 0\}$ and satisfying \eqref{Edomain} is $D_t=B(\tilde 0,R_t)$ with $(R_t)_t$ a diffusion process starting at $0$ which infinitesimal generator is $L_n:=\partial_r^2+b_n(r)\partial_r$ where $b_n(r)=\frac{d}{dr}\ln\left(\frac{I_n(r)^2}{I_n'(r)}\right)$ and $I_n$ is given by \begin{equation}\label{EIn}
    I_n(r)=\begin{cases}
    \int_0^r\sin^{n-1}(s)ds & \text{ if } M_n=\mathbb{S}^n\\
     \int_0^r\sin^{n-1}\left(\frac{s}{2}\right)ds & \text{ if } M_n=\mathbb{RP}^n\\
     \int_0^{r}  \sin^{2n-1}\left(\frac{s}{2}\right)\cos\left(\frac{s}{2}\right) \, ds=
 \frac{ \sin^{2n}\left(\frac{r}{2}\right)}{n} &\text{ if }M_n=\mathbb{CP}^n\\
  \int_0^{r}  \sin^{4n-1} \left(\frac{s}{2}\right)\cos^3 (\frac{s}{2})\, ds= \frac{ \sin^{4n} \left(\frac{r}{2}\right)( 1 + 2n \cos^2 (\frac{r}{2}))}{2n(2n+1)} &\text{ if }M_n=\mathbb{HP}^n
    \end{cases}.\end{equation}
    Moreover, the first covering time $\tau_n:=\inf\{t>0|\, R_t=\pi\}$ of $(D_t)_{t \in [0, \tau_n]}$ is a.s. finite.
\end{pro}
\proof
From Proposition \ref{PVHarmonic}, \eqref{Ederive} can be rewritten $b_n(r)=\frac{d}{dr}\ln\left(\frac{I_n(r)^2}{I_n'(r)}\right)$. Using the notations from Section 6 in~\cite{Karlin}, to show that $0$ is an entrance boundary, it is enough to show that 
    \begin{equation}\label{EntranceBound}
        S(0,\pi/2):=\int_0^{\pi/2}\frac{I_n'(r)}{I_n^2(r)}dr=+\infty
        \text{ and }N(0,\pi/2):=\int_0^{\pi/2}\frac{I_n^2(r)}{I_n'(r)}\int_r^{\pi/2}\frac{I_n'(s)}{I_n^2(s)}ds\,dr<+\infty.
    \end{equation}
    Using the fact that $I_n$ is non decreasing, we can remark that 
    $$N(0,\pi/2)=\int_0^{\pi/2} \frac{I_n(r)}{I_n'(r)}\frac{I_n(\pi/2)-I_n(r)}{I_n(\pi/2)}dr\leq \int_0^{\pi/2} \frac{I_n(r)}{I_n'(r)}dr.$$
 From \eqref{EIn}, for $\mathbb{CP}^n$ and $\mathbb{HP}^n$ we directly get
 $\frac{I_n'(r)}{I_n^2(r)}\underset{r\to 0}{\sim}\left(\frac{an}{2}\right)^2\left(\frac{2}{r}\right)^{an+1}$ and $\frac{I_n(r)}{I_n'(r)}\underset{r\to 0}{\sim}\frac{r}{an}$
which is enough to show \eqref{EntranceBound}.
For the real projective space it is enough to remark that $\int_0^r\sin^{n-1}\left(\frac{s}{2}\right)ds\leq r\sin^{n-1}\left(\frac{r}{2}\right)$ to obtain $\frac{I_n'(r)}{I_n^2(r)}\geq \frac{1}{r^2\sin^{n-1}\left(\frac{r}{2}\right)}$ and $\frac{I_n(r)}{I_n'(r)}\leq r$ and to conclude. The sphere case works in the same way.
The last remark about the first covering time is postponed to Remark \ref{Rem1}.
\wwtbp

\begin{Rem}
Since $ \tau_n$ is a.s. finite, we extend $ D $ after time $ \tau_n$, by $D_s= D_{\tau_n} = M_n$ for $\tau_n \le s $ and consider $ (D_t)_{t \ge 0} $. 
\end{Rem}

 Let $X^n $ be a Brownian motion in $M_n $ started at $\tilde 0 \in M_n $, and accelerated by a factor $2$.
 In  \cite{ arnaudon:hal-03037469}, several couplings of  $X^n$ and $D$ were constructed, so that for any time $t\geq 0$, the conditional law of $X^n_t$ knowing the trajectory $D_{[0,t]}\df(D_s)_{s\in[0,t]}$ is the normalized uniform law over $D_t$, which is $\Lambda(D_t,\cdot)$, i.e.
\bqn{loiprod}
\fo t\geq 0,\qquad \cL(X_t\vert  D_{[0,t]})&=&\Lambda(D_t,\cdot),\eqn
\par
 where $ \Lambda$ is the corresponding Markov kernel defined in the beginning of Section~\ref{Sec0}.
 
 In the following  corollary we explicit a coupling that satisfies \eqref{loiprod}, which were constructed in \cite{arnaudon:hal-03037469} Theorems 3.5. Before stating it, we fix some notations. Recall that the cut-locus of $\tilde 0$ is given by $\text{cut}(\wi 0) = \{ x \in M_n,\text{ s.t. }\, d(\wi 0, x)= \pi \}$. For $x\in M_n \backslash\{\wi 0, \text{cut}(\wi 0)\}$ we denote by $N^{\scriptscriptstyle in}(x)$ the unit vector at $x$ normal to the sphere centred at $\wi 0$ with radius $d(\wi 0,x)$ (where $d $ is the distance in $M_{n}$) and pointing towards $\wi 0$, i.e., $N^{\scriptscriptstyle in}(x)=-\nabla d(\wi 0,\cdot)(x)$.\\ 
\begin{cor}\label{cor1}
Consider the Brownian motion $X^n\df (X^n_t)_{t\geq 0}$ in $M_{n}$  accelerated by a factor $2$.\\
 \textbf{Full coupling}.
 Let $D^{(1)}_t$ be the ball in $M_{n}$ centred at $\wi 0$ with radius $R^{(1)}_t$ solution started at $0$ to the It\^o equation
\bq
dR^{(1)}_t&=&- \langle N {^{\scriptscriptstyle in}}(X^n_t), dX^n_t)\rangle + \left[2 \frac{I_{n}''}{I_{n}'}(d(\wi 0,X^n_t))-\frac{I_{n}''}{I_{n}'}(R^{(1)}_t)\right]\, dt
\eq
This evolution equation is considered up to  the hitting time  $\tau_n^{(1)}$ of $\pi$ by $R^{(1)}_t$. \\

Let $D_t$ be the ball  in $M_{n}$ centered at $\wi 0$ with radius $R_t$ defined in Proposition \ref{T01}, and let $\tau_n$ be the stopping time defined in Proposition \ref{T01}.

Then we have: 
\begin{itemize}
\item[(1)]
The coupling $(X_t,D_t^{(1)})_ {0\leq t\leq \tau_n^{(1)}}$ is satisfying \eqref{loiprod}, changing $t$ by any $ \mathcal{F}^{D^{(1)}}$-stopping time $ \tau \le \tau_{n}^{(1)}$;
\item[(2)] the pairs  $\left(\tau_n^{(1)}, D^{(1)}_{ [0,\tau_n^{(1)}]}\right)$, and  $\left(\tau_n,D_{ [0,\tau_n]}\right)$ have the same law.
\end{itemize}
\end{cor}
 
\begin{Rem}
Note that   $ \tau_n^{(1)}$ is a stopping time for $ \mathcal{F}^{D^{(1)}} $ and so for $ \mathcal{F}^{X}$.
\end{Rem} 

 Due to these couplings and to general arguments from Diaconis and Fill \cite{MR1071805} for discrete case (for the continuous case see Proposition
 \ref{sep} below), $\tau_{n}^{(1)}$ is a strong stationary time for $X^n$, meaning that $\tau_{n}^{(1)}$ and $X^n_{ \tau_{n}^{(1)}}$ are independent and $X^n_{ \tau_{n}^{(1)}}$ is uniformly distributed over  $M_{n}$.

\begin{pro}\label{sep} If $\tau_{n}$ and so $\tau_{n}^{(1)} $ is finite almost surely, then $\tau_{n}^{(1)}$ is a strong stationary time for $X^n$, and 
\bq
\fo t\geq 0,\qquad \fs(\cL(X^n_t),\cU_{n})&\leq & \PP[\tau_{n}^{(1)} > t] = \PP[\tau_{n} > t].
\eq
\end{pro}
\proof
For simplicity write $X_t := X^n_t $. 
Let $ f : M_n \to \mathbb{R} $ be a bounded measurable function,  apply \eqref{loiprod} at $ \tau_{n}^{(1)}$, since $D^{(1)}_{\tau_{n}^{(1)}} = B(\wi 0, \pi) = M_n$ we get :
\begin{align*}
\EE[ f(X_{\tau_{n}^{(1)} })] &=  \EE[  \EE[ f(X_{ \tau_{n}^{(1)}}) \vert  D^{(1)}_{[0,\tau_{n}^{(1)}]} ] ]\\
&=   \EE[  \Lambda(D^{(1)}_{\tau_{n}^{(1)} }, f)]\\
&= \EE[ \Lambda (M_n, f) ] \\
&=  \frac{\int_{M_n} f \, d V_g }{V_g(M_n)}   = \cU_{n}(f) .
\end{align*}
 Hence $X_{\tau_{n}^{(1)}}$ is uniformly distributed over  $M_{n}$, and clearly $\tau_{n}^{(1)} $ and $X_{\tau_{n}^{(1)}}$ are independent so $\tau_{n}^{(1)}$ is a strong stationary time for $X^n$. \\
Let $ f : M^n \to \mathbb{R}_+ $ be a bounded positive measurable function,   and  $  0 < t$,  we have :
 \begin{align}
\EE_{\wi 0}\left[ f(X_t)\right] & \ge  \EE_{(\wi 0,0)}\left[ f(X_t) \un_{\tau_{n}^{(1)} \le t} \right] \text{ where }(\tilde 0,0)\text{ is the starting point for }(X,D^{(1)})\notag\\ 
 & =  \EE_{(\wi 0,0)}\left[ \un_{ \tau_{n}^{(1)} \le t} \EE \left[f(X_t) \vert \mathcal{F}_ {\tau_{n}^{(1)}}\right] \right] \notag \\ 
 &=  \EE_{(\wi 0,0)}\left[ \un_{ \tau_{n}^{(1)} \le t} \EE_{X_{\tau_{n}^{(1)}} } \left[f(X_{t- \tau_{n}^{(1)} })\right] \right]  \quad \text{Strong Markov property}\\ 
 &= \cU_{n}(f) \mathbb{P}_{0} \left[\tau_{n}^{(1)} \le t \right],\label{Esep} 
\end{align}
where in the last line we use that $ X_{\tau_{n}^{(1)} } \sim\cU_{n}$   is invariant under  $X$ and that $  X_{ \tau_{n}^{(1)}}$ and $\tau_{n}^{(1)} $  are independent.

Let $y \in M_n$, $\epsilon >0$. Let  $f= \un_{B(y, \epsilon)}   $, and let $p_t( \wi 0 ,x) $ be the heat kernel, i.e., $X_t(\wi 0) \overunderset{\mathcal{L}}{}{\sim} p_t( \wi 0 ,x) \,\cU_{n} (dx)$,  then \eqref{Esep} becomes
$$ \frac{  \int  p_t( \wi 0 ,x) \un_{B(y, \epsilon)}  (x) \cU_{n} (dx) }{ \cU_{n}  (B(y, \epsilon))} \ge \mathbb{P}_{0} [  \tau_{n}^{(1)} \le t ], $$ 
letting $ \epsilon $ goes to $0$, we obtain for all $y \in M_n$,
$$ p_t( \wi 0 ,y)  \ge \mathbb{P}_{0} [\tau_{n}^{(1)} \le t ],$$
hence 
$$ 1 - p_t( \wi 0 ,y) \le  \mathbb{P}_{0} [\tau_{n}^{(1)}  > t ].$$
and by definition of the separation discrepancy the result follows.

\wwtbp   \vspace{0.5cm}

\par
\begin{Rem}
Note that for any $t\in [0,\tau_{n}^{(1)} )$, the ``opposite pole'' $cut(\wi 0)$ does not belong to the support of $\Lambda(D^{(1)}_t,\cdot)$.
It follows from an extension of Remark 2.39 of Diaconis and Fill \cite{MR1071805} that $ \tau_{n}^{(1)} $ is even a sharp  strong stationary time for $X^n$, meaning
that
\bqn{sharps}
\fo t\geq 0,\qquad \fs(\cL(X^n_t),\cU_{n})&= & \PP[ \tau_{n}^{(1)} > t].\eqn
\par
Thus for all $M_n$ the understanding of the convergence in separation of $X^n$ toward $\cU_{n}$ amounts to understanding the tail  distribution of $\tau_{n}$. 
\end{Rem}

\subsection{Study of the first covering time}
To study separation cutoff phenomenon we want to obtain good estimates for the expectation and the variance of the first covering time  $\tau_n$. As in~\cite{ACM}, this can be done by considering the Green operator $G_n$ associated to the Poisson equation:
\begin{equation}
    \begin{cases}
    &L_n\phi_n=-g,\,g\in\mathcal{C}_b([0,\pi])\\
    &\phi_n(\pi)=0, \quad \phi_n \text{ is bounded}
    \end{cases}
\end{equation} where $L_n$ still denotes the infinitesimal generator defined in Proposition \ref{T01}.
The Green operator is given by:
\begin{equation}
\begin{array}{ccccc}
G_n & : & \mathcal{C}_b([0,\pi]) & \to & \mathcal{C}_b([0,\pi]) \\
 & & g & \mapsto & \int_{\cdot}^{\pi}\frac{I_n'(t)}{I_n^2(t)}\int_0^t\frac{I_n^2(s)}{I_n'(s)}g(s)dsdt \\
\end{array}
\end{equation}
Note that the definition of $G_n$ is ensured using the same arguments as for the finiteness of $N(0,\pi/2)$ in~\eqref{EntranceBound}.

As for Propositions 7 and 9 in~\cite{ACM} the moments and the variance of $\tau_n$ can then be expressed in function of this Green operator.
\begin{pro}
For $k\geq 0$, set $u_{n,k}=k!G_n^{\circ k}[\mathbb{1}]$ where $G_n^{\circ k}$ denotes the $k$-th iteration of $G_n$. Then, for all $n\geq 1$ and $k\geq 0$: \begin{equation}\label{Emoments}
\EE[\tau_n^k]=u_{n,k}(0).
\end{equation}
In particular
\begin{equation}\label{EEsp}
    \EE[\tau_n]=\int_0^\pi\frac{I_n(t)}{I_n'(t)}\frac{I_n(\pi)-I_n(t)}{I_n(\pi)}dt
\end{equation}
\begin{align}
\Var(\tau_n)&=2\int_0^{\pi}\frac{I_n'(t)}{I_n(t)^2}\int_0^t\frac{I_n^2(s)}{I_n'(s)}u_{n,1}'(s)^2dsdt\label{EVar}\\
&=\frac{2}{I_n(\pi)^2} \int_0^{\pi}\frac{I_n^2(u)}{I'_n(u)} \int_u^{\pi}\frac{(I_n(\pi)-I_n(s))^2}{I'_n(s)}  ds \, du.\label{E0}
\end{align}

\end{pro}
\proof
The proof for \eqref{Emoments}, \eqref{EEsp} and \eqref{EVar} is the same than in~\cite{ACM}. For \eqref{E0}, we remark that

\[
\left( u'_{n,1}(s) \right)^2 = \left( \frac{I'_n(s)}{I_n^2(s)} \int_0^s \frac{I_n^2(u)}{I'_n(u)} du \right)^2
\]
thus
\begin{align*}
\int_0^t \frac{I_n^2(s)}{I_n'(s)} \left( u'_{n,1}(s) \right)^2 ds
&=\int_0^t\frac{I_n'(s)}{I_n^2(s)}\left(\int_0^s\frac{I_n^2(u)}{I_n'(u)}du\right)^2ds\\ 
&= \left[ \frac{-1}{I_n(s)} \left(\int_0^s \frac{I_n^2(u)}{I'_n(u)} du  \right)^2 \right]_0^t
+ 2 \int_0^t \frac{I_n(s)}{I'_n(s)} \int_0^s \frac{I_n^2(u)}{I'_n(u)} du ds\\
&= \frac{-2}{I_n(t)} \int_0^t \frac{I_n^2(s)}{I'_n(s)} \int_0^s \frac{I_n^2(u)}{I'_n(u)} du ds
+ 2 \int_0^t \frac{I_n(s)}{I'_n(s)}\int_0^s \frac{I_n^2(u)}{I'_n(u)} du ds\\
&= 2 \int_0^t \frac{I_n(s)}{I'_n(s)} \left( 1 - \frac{I_n(s)}{I_n(t)} \right)
 \int_0^s \frac{I_n^2(u)}{I'_n(u)} du \, ds
\end{align*}

Using Fubini's theorem and integrating by parts:
\begin{align*}
\mathrm{Var}(\tau_n) 
&=4 \int_0^{\pi} \frac{I'_n(t)}{I_n^2(t)}
\int_0^t \frac{I_n(s)}{I'_n(s)} \left(1 - \frac{I_n(s)}{I_n(t)} \right)
\int_0^s \frac{I_n^2(u)}{I'_n(u)} du  ds \, dt
\\
&= 4 \int_0^{\pi}\frac{I_n^2(u)}{I'_n(u)} \int_u^{\pi}\frac{I_n(s)}{I'_n(s)} \left(\int_s^{\pi}\frac{I'_n(t)}{I_n^2(t)} - I_n(s)\int_s^{\pi}\frac{I'_n(t)}{I_n^3(t)}\right) dt  ds \, du\\
&= 4 \int_0^{\pi}\frac{I_n^2(u)}{I'_n(u)} \int_u^{\pi}\frac{I_n(s)}{I'_n(s)} \left(\left[\frac{I_n(s)-2I_n(t)}{2I_n(t)^2}\right]_s^{\pi}\right)  ds \, du\\
&= 4 \int_0^{\pi}\frac{I_n^2(u)}{I'_n(u)} \int_u^{\pi}\frac{I_n(s)}{I'_n(s)} \left(\frac{I_n(s)-2I_n(\pi)}{2I_n(\pi)^2}+\frac{1}{2I_n(s)}\right) ds \, du\\
&= \frac{2}{I_n(\pi)^2} \int_0^{\pi}\frac{I_n^2(u)}{I'_n(u)} \int_u^{\pi}\frac{(I_n(\pi)-I_n(s))^2}{I'_n(s)}  ds \, du.
\end{align*}
\wwtbp

\begin{Rem}\label{Rem1}
In particular this proves that $\EE[\tau_n]<+\infty$ and thus $\tau_n<+\infty$ a.s.
\end{Rem}
\section{ Cutoff for the complex projective space}\label{Sec1}
\begin{pro}\label{PR1}
On the complex projective space, the first covering time $\tau_n$ is satisfying:
\begin{align*}
    \mathbb{E}[\tau_n]& =\frac{1}{n}\sum_{k=1}^{n} \frac{1}{k}= \frac{\ln(n)}{n} + \frac{\gamma}{n}+ \frac{1}{2n^2} + o\left(\frac{1}{n^2}\right)\\
    \Var(\tau_n)&=\frac{\pi^2}{6n^2}+o\left(\frac{1}{n^2}\right)
\end{align*}
where $\gamma$ is the Euler–Mascheroni constant.
\end{pro}
This proposition will be proven in the two next subsections (Subsections \ref{Subsec1} and \ref{Subsec2}).

From this and using \eqref{sharps}, we get the following estimate for the cutoff :

\begin{theo}[For the complex projective space]\label{T1} Let $X^n$ be the (two-times accelerated) Brownian Motion in $ \mathbb{CP}^n$ and $\gamma$ the Euler–Mascheroni constant.
Let $ c \in \mathbb{R} $ and $t_n = \frac{\ln(n)}{n} + \frac{c}{n}$:
\begin{itemize}
\item if  $ c > \gamma$, 
  $\limsup\limits_{n \to \infty}  \fs(\cL(X^n_{t_n}),\mathcal{U}_{n}) \le \frac{\pi^2}{6 (c- \gamma)^2}   $
\item  if  $ c < \gamma$ then
  $\liminf\limits_{n \to \infty}  \fs(\cL(X^n_{t_n}),\mathcal{U}_{n})  \ge 1- \frac{\pi^2}{6 (c- \gamma)^2}   $.
\end{itemize}
Hence  $X^n$ has a cutoff in separation at time $\frac{\ln(n)}{n}  $ with  window $ \frac1n$, and we have the above control.
\end{theo}
\proof

For $ c > \gamma$, for $n$ large enough, $t_n-\EE[\tau_n]>0$, thus:
\bq
\fs(\cL(X^n_{t_n}),\mathcal{U}_{n}) &= & \PP[\tau_n>t_n]  \\
                 &= & \PP[\tau_n - \mathbb{E}[ \tau_n ] >t_n - \mathbb{E}[ \tau_n ] ] \\
                 & \le & \frac{\mathrm{Var}(\tau_n)} {(t_n - \mathbb{E}[ \tau_n])^2} \\
                 &=& \frac{\frac{\pi^2}{6n^2} + o\left(\frac{1}{n^2}\right)}{ \left( \frac{c- \gamma}{n}+o(\frac{1}{n}) 
                 \right)^2}  \\
                 &=& \frac{\frac{\pi^2}{6} + o\left(1 \right)}{ \left( {c- \gamma}+ o(1) \right)^2}.
\eq

Hence, for $ c > \gamma$,
  $$\limsup_{n \to \infty}  \fs(\cL(X^n_{t_n}),\mathcal{U}_{n}) \le \frac{\pi^2}{6 (c- \gamma)^2} .  $$

For $ c < \gamma$, $t_n-\EE[\tau_n]<0$ for $n$ large enough, thus:
\bq
\fs(\cL(X^n_{t_n}),\mathcal{U}_{n})&= & 1 - \PP[\tau_n \le t_n]  \\
                 & = & 1 - \PP[\tau_n - \mathbb{E}[ \tau_n ] \le t_n - \mathbb{E}[ \tau_n ] ] \\ 
                & = & 1 - \PP\left[\tau_n - \mathbb{E}[ \tau_n ] \le \frac{c- \gamma}{n}+ o\left(\frac1n\right) \right] \\ 
                 &\ge & 1- \frac{\frac{\pi^2}{6n^2} + o\left(\frac{1}{n^2}\right)}{ \left( \frac{c- \gamma}{n}+ o(\frac{1}{n})
                 \right)^2}  \\
                 &=& 1- \frac{\frac{\pi^2}{6} + o\left(1 \right)}{ \left( {c- \gamma}+ o(1) \right)^2}.
                 \eq
Hence for $ c < \gamma$,
  $$\liminf_{n \to \infty}  \fs(\cL(X^n_{t_n}), \mathcal{U}_{n}) \ge 1- \frac{\pi^2}{6 (c- \gamma)^2}   $$
  
\wwtbp
\subsection{Computation of the mean, first part of the proof of Proposition \ref{PR1}}\label{Subsec1}

For the complex projective space ($a=2$), from \eqref{EIn} we remind that :

\[
I_n(r)  
= \frac{ \sin^{2n}\left(\frac{r}{2}\right)}{n},\,I_n'(r)=\cos\left(\frac{r}{2}\right)\sin^{2n-1}\left(\frac{r}{2}\right)
\quad \text{and} \quad
I_n(\pi) = \frac{1}{n}.
\]
Then
\begin{align*}
\mathbb{E}[\tau_n] &= \int_0^{\pi} \frac{\sin^{2n}\left(\frac{x}{2}\right)}{n \, \sin^{2n-1}\left(\frac{x}{2}\right) \cos\left(\frac{x}{2}\right)} \left(1 - \sin^{2n}\left(\frac{x}{2}\right)\right) dx\\
&= \frac{2}{n} \int_0^{\pi/2} \frac{\sin(x)}{\cos(x)} \left(1 - \sin^{2n}(x)\right) dx\\
&= \frac{2}{n} \int_0^1 \frac{y}{1 - y^2} \left(1 - y^{2n} \right) dy
= \frac{2}{n} \int_0^1 \sum_{k=0}^{n-1} y^{2k+1} \, dy
= \frac{1}{n} \sum_{k=1}^{n} \frac{1}{k}.
\end{align*}

Recall that
\begin{equation}\label{Elog}
\sum_{k=1}^{n} \frac{1}{k} = \ln(n) + \gamma + \frac{1}{2n} + o\left(\frac{1}{n}\right)
\end{equation}
where $\gamma$ is the Euler–Mascheroni constant.

We deduce that :
\[
\mathbb{E}[\tau_n] = \frac{\ln(n)}{n} + \frac{\gamma}{n}+ \frac{1}{2n^2} + o\left(\frac{1}{n^2}\right).
\]

\subsection{The computation of the variance, second part of the proof of Proposition \ref{PR1}}\label{Subsec2}

From \eqref{E0}, we have
\begin{equation*}
\mathrm{Var}(\tau_n) 
=\frac{8}{n^2} \int_0^{\frac{\pi}{2}}\frac{\sin^{2n+1}(u)}{\cos\left(u\right)} \int_u^{\frac{\pi}{2}}\frac{\left(1-\sin^{2n}(s)\right)^2}{\sin^{2n-1}(s)\cos(s)}  ds \, du.\\
\end{equation*}
After the following change of variables:
\[X=\sin(u),\ Y=\sin(s)\]
we get :
\begin{equation*} \mathrm{Var}(\tau_n)=\frac{8}{n^2} \int_0^1  \frac{X^{2n+1}}{1 - X^2} \int_X^1 \frac{Y}{1 - Y^2} \left( \frac{1 - Y^{2n}}{Y^{n}} \right)^2 \, dY \, dX
\end{equation*}

Also
\begin{align*}
\int_X^1 \frac{Y}{1 - Y^2} &\left( \frac{1 - Y^{2n}}{Y^{n}} \right)^2 dY \\
&= \int_X^1 \frac{1}{Y^{2n - 1}} (1 - Y^{2n})\sum_{k=0}^{n-1} Y^{2k} \, dY
\\
&= \sum_{k=0}^{n-1} \int_X^1 Y^{2(k - n) + 1} - Y^{2k + 1} \, dY
\\
&= \sum_{k=0}^{n-2} \frac{1}{2(k - n +1)} \left( 1 - X^{2(k - n +1)} \right) - \sum_{k=1}^{n} \frac{1}{2k} (1 - X^{2k}) - \ln(X)
\\
&= -\frac{1}{2} \sum_{k=1}^{n-1} \frac{1}{k} \left( (1 - X^{-2k}) + (1 - X^{2k}) \right) - \frac{1}{2n} (1 - X^{2n}) - \ln(X)
\\
&= \frac{1}{2} \sum_{k=1}^{n-1} \frac{1}{k} \left(  \frac{1-X^{2k}}{X^k}  \right)^2 - \frac{1}{2n}(1 - X^{2n}) - \ln(X)
\end{align*}

\begin{align*}
\mathrm{Var}(\tau_n) &= \frac{8}{n^2} \int_0^1  \frac{X^{2n+1}}{1 - X^2} 
\left[ \frac{1}{2} \sum_{k=1}^{n-1} \frac{1}{k} \left( \frac{1 - X^{2k}}{X^k} \right)^2 
- \frac{1}{2n}(1 - X^{2n}) - \ln(X) \right] dX
\\
&=\underbrace{\frac{4}{n^2} \sum_{k=1}^{n-1} \frac{1}{k} \int_0^1  \frac{X^{2n+1}}{1 - X^2} \left( \frac{1 - X^{2k}}{X^k} \right)^2 dX }_{A_n}\\
&- \underbrace{\frac{4}{n^3} \int_0^1  \frac{X^{2n+1}}{1 - X^2}(1 - X^{2n}) dX}_{B_n} \\
& - \underbrace{ \frac{8}{n^2} \int_0^1  \frac{X^{2n+1}}{1 - X^2}\ln(X)  dX}_{C_n}.
\end{align*}

Let us compute an asymptotic expansion of $A_n$:
\begin{align*}
A_n &= \frac{4}{n^2} \sum_{k=1}^{n-1} \frac{1}{k} 
\int_0^1 \frac{X^{2n+1 - 2k}(1 - X^{2k})^2}{1 - X^2} \, dX\\
&= \frac{4}{n^2} \sum_{k=1}^{n-1} \frac{1}{k} \sum_{l=0}^{k-1}
\int_0^1 X^{2n+1 - 2k}(1 - X^{2k}) X^{2l}\, dX\\
&= \frac{4}{n^2} \sum_{k=1}^{n-1} \frac{1}{k} \sum_{l=0}^{k-1}
\int_0^1 X^{2n+1 - 2k+2l} -  X^{2l+2n+1}\, dX\\
&= \frac{2}{n^2} \sum_{k=1}^{n-1} \frac{1}{k} 
\sum_{l=0}^{k-1} \left( \frac{1}{n + l - k + 1} - \frac{1}{n + l + 1} \right)\\
&= \frac{2}{n^2} \sum_{k=1}^{n-1} \sum_{l=0}^{k-1} 
\frac{1}{k} \left( \frac{1}{n + l - k + 1} - \frac{1}{n + l + 1} \right)\\
&= \frac{2}{n^2} \sum_{k=1}^{n-1} \sum_{l=0}^{k-1} 
\left( \frac{1}{(n + l - k + 1)(n + l + 1)} \right)\\
&=\frac{2}{n^2} \sum_{k=1}^{n-1} \sum_{l=1}^{k} 
\left( \frac{1}{\left(1 + \frac l n - \frac k n \right)\left(1 + \frac l n \right)} \frac{1}{n^2}\right)
\end{align*}
Now notice that as a Riemann sum:
\[
\sum_{k=1}^{n-1} \sum_{l=1}^{k} \frac{1}{\left(1 + \frac{l}{n} - \frac{k}{n}\right)(1 + \frac{l}{n})}  \frac{1}{n^2}
\longrightarrow \int_0^1 \int_0^x \frac{1}{(1 + y - x)(1 + y)} \, dy \, dx \quad \text{as } n \to \infty
\]

and
\begin{align*}
\int_0^1 \int_0^x \frac{1}{(1 + y - x)(1 + y)} \, dy \, dx 
&= \int_0^1 \left( \int_y^1 \frac{1}{(1 + y - x)} \, dx \right) \frac{1}{1 + y} \, dy
= - \int_0^1 \frac{\ln( y)}{1+y} \, dy\\
&= - \sum_{k=0}^{\infty} (-1)^k \int_0^1 \ln( y)y^k \,dy
= \sum_{k=0}^{\infty} \frac{(-1)^k}{(k+1)^2}
= \frac{\pi^2}{12}
\end{align*}
where we used integration by parts (for all $k \ge 0$) to obtain
\begin{equation}\label{logIPP}
\int_0^1 y^k \ln(y) \, dy = -\frac{1}{(k+1)^2}.
\end{equation}
Hence, we conclude that:
\[
A_n \sim \frac{\pi^2}{6n^2} \quad \text{as } n \to \infty.
\]

Let us finally compute an asymptotic expansion of $B_n$:

\begin{align*}
B_n &= \frac{4}{n^3} \int_0^1 \frac{X^{2n+1}}{1 - X^2} \left(1 - X^{2n} \right) \, dX\\
&= \frac{4}{n^3} \sum_{k=0}^{n-1} \int_0^1 X^{2n + 2k + 1} \, dX\\
&= \frac{2}{n^3} \sum_{k=0}^{n-1} \frac{1}{n + k + 1}
= \frac{2}{n^3} \sum_{k=n+1}^{2n} \frac{1}{k}.
\end{align*}

Hence, using the Euler–Mascheroni asymptotic we get 
\[
B_n \sim \frac{2}{n^3} \ln(2).
\]

Let us finally compute the asymptotic expansion of $C_n$:

\[
C_n = \frac{8}{n^2} \int_0^1 \frac{X^{2n+1}}{1 - X^2} \ln(X) \, dX
= \frac{8}{n^2} \sum_{k=0}^{\infty} \int_0^1 X^{2n+1 + 2k} \ln(X) \, dX.
\]

Using again \eqref{logIPP}:

\[
C_n = -\frac{8}{n^2} \sum_{k=0}^{\infty} \frac{1}{(2n + 2k + 2)^2}
= -\frac{2}{n^2} \sum_{k=n+1}^{\infty} \frac{1}{k^2}
\sim -\frac{2}{n^3}
\]

and we also have:
\[
\mathrm{Var}(\tau_n) = \frac{\pi^2}{6n^2} + o\left(\frac{1}{n^2}\right).
\]

\section{Cutoff for the  quaternionic projective space}\label{Sec2}
\begin{pro}\label{PR2}
On the quaternionic projective space, the first covering time $\tau_n$ is satisfying:
\begin{align*}
    \mathbb{E}[\tau_n]& = \frac{1}{2n+1} \sum_{k=1}^{2n+1} \frac{1}{k} = \frac{\ln(2n+1)}{2n+1} + \frac{\gamma}{2n+1}+ \frac{1}{2(2n+1)^2} + o\left(\frac{1}{n^2}\right)\\
    \Var(\tau_n)&\underset{n\to+\infty}{\sim}\frac{N}{(2n+1)^2}
\end{align*}
where $\gamma$ is the Euler–Mascheroni constant and 
\begin{equation*}
    N=2\int_{0\leq y\leq x<\infty} e^{y-x}\frac{\left(x+1\right)^2}{x^2}\frac{\left(1-e^{-y}(y+1)\right)^2}{y^2}dx\,dy.
\end{equation*}
\end{pro}
As previously, the proof will be postponed to the two following subsections (Subsections \ref{Subsec3} and \ref{Subsec4}). With the same method than for Theorem \ref{T1} we obtain the following result for the cutoff:
\begin{theo}[For the quaternionic projective space]\label{T2} Let $X^n$ be the (two-times accelerated) Brownian motion in $ \mathbb{HP}^n$, $\gamma$ the Euler–Mascheroni constant and $N$ the constant obtained in Proposition \ref{PR2}.
Let $ c \in \mathbb{R} $ and $t_n = \frac{\ln(2n+1)}{2n+1} + \frac{c}{2n+1}$: 
\begin{itemize}
\item if  $ c > \gamma$ then
  $\limsup\limits_{n \to \infty}  \fs(\cL(X^n_{t_n}),\mathcal{U}_{n}) \le \frac{N}{(c- \gamma)^2}   $
\item  if  $ c < \gamma$ then
  $\liminf\limits_{n \to \infty}  \fs(\cL(X^n_{t_n}),\mathcal{U}_{n}) \ge 1- \frac{N}{ (c- \gamma)^2}   $.
\end{itemize}
Hence  $X^n$ has a cutoff in separation at time $\frac{\ln(2n+1)}{2n+1}  $ with  window $ \frac{1}{2n+1}$, and we have the above control.
\end{theo}

\subsection{Computation of the mean, first part of the proof of Proposition \ref{PR2}}\label{Subsec3}

For the quaternionic projective space ($a=4$), \eqref{EIn} gives:

\begin{equation}\label{E4}
\begin{split}I_n(r)&
= \frac{ \sin^{4n} \left(\frac{r}{2}\right)( 1 + 2n \cos^2 (\frac{r}{2}))}{2n(2n+1)}\quad \text{and} \quad
I_n(\pi) = \frac{1}{2n(2n+1)}\\ I_n'(r)&{=\cos^3\left(\frac{r}{2}\right)\sin^{4n-1}\left(\frac{r}{2}\right)}.\end{split}
\end{equation}

\begin{align*}
\mathbb{E}[\tau_n] &= \frac{1}{2n(2n+1)} \int_0^{\pi} \frac{\sin\left(\frac{x}{2}\right)}{ \cos^3\left(\frac{x}{2}\right)} \left(1 + 2n \cos^{2} \left(\frac{x}{2}\right)\right)\left( 1 - \sin^{4n} \left(\frac{x}{2}\right)( 1 + 2n \cos^2\left(\frac{x}{2}\right) \right)dx \\
&=\frac{1}{2n(2n+1)} \int_0^{1} \frac{P(X)} {(1-X)^2} dX=\frac{1}{2n(2n+1)}\left(\left[\frac{P(X)}{1-X}\right]_0^1-\int_0^1\frac{P'(X)}{1-X}dX\right)
\end{align*}
with the change of variable $X=\sin^2\left(\frac x 2\right)$ and 
$$ P(X) = {\left(1 + 2n (1 -  X) \right)}  \left( 1 - X^{2n}  (1 + 2n (1 -  X)) \right).$$
Note that $\left[\frac{P(X)}{1-X}\right]_0^1=-(2n+1)$.
Since
\begin{equation*}
P'(X) = -2n \left( 1 - X^{2n}  (1 + 2n(1 -  X)) \right) -2n ( 2n+1) (1 + 2n (1 -  X)) X^{2n-1}(1-X  ):
\end{equation*}
\begin{align*}
- \int_0^{1} \frac{ P'(X)}{1-X} dX  &= 2n \int_0^{1} \left( \frac{1 - X^{2n}}{1-X}     -2nX^{2n} + (2n+1)(1 + 2n (1 -  X)) X^{2n-1} \right) dX\\
&= 2n\left(    \sum_{k=1}^{2n} \frac{1}{k} - \frac{2n}{2n+1} + \frac{2n+1}{2n} + 2n(2n+1) ( \frac{1}{2n} - \frac{1}{2n+1})  \right)\\
&= 2n\left(    \sum_{k=1}^{2n} \frac{1}{k} - \frac{2n}{2n+1} + \frac{2n+1}{2n} + 1  \right)\\
\end{align*}
Finally
\begin{align*}
\mathbb{E}[\tau_n] &= \frac{1}{2n(2n+1)} \left( -(2n+1)     +    2n\left(    \sum_{k=1}^{2n} \frac{1}{k} - \frac{2n}{2n+1} + \frac{2n+1}{2n} + 1  \right)      \right) \\
&= \frac{1}{2n(2n+1)} \left(     2n\left(    \sum_{k=1}^{2n} \frac{1}{k} - \frac{2n}{2n+1} + 1  \right)      \right) \\
&= \frac{1}{2n+1} \sum_{k=1}^{2n+1} \frac{1}{k}    \\
&= \frac{\ln(2n+1)}{2n+1} + \frac{\gamma}{2n+1}+ \frac{1}{2(2n+1)^2} + o\left(\frac{1}{n^2}\right).
\end{align*}

\subsection{The computation of the variance, second part of the proof of Proposition \ref{PR2}}\label{Subsec4}
Applying \eqref{E4} to \eqref{E0}, we have 
\begin{equation*}
  \begin{split}
     \Var(\tau_n)=\frac{8}{\left(2n(2n+1)\right)^2}\int_0^{\frac{\pi}{2}}&\frac{\sin^{4n+1}(u)\left(2n\cos^2(u)+1\right)^2}{\cos^3(u)}\\&\times\int_u^{\frac{\pi}{2}}\frac{\left(1-\sin^{4n}(s)\left(2n\cos^2(s)+1\right)\right)^2}{\sin^{4n-1}(s)\cos^3(s)}ds\,du. 
  \end{split}  
\end{equation*}
With the change of variables $x=2n\cos^2(u)$ and $y=2n\cos^2(s)$:
\begin{equation*}
    \Var(\tau_n)=\frac{2}{\left(2n+1\right)^2}\int\limits_{[0,+\infty[^2}F_n(x,y)dx\,dy.
\end{equation*}
with \begin{equation*}
    F_n(x,y):=\frac{\left(1-\frac{x}{2n}\right)^{2n}\left(x+1\right)^2}{x^2}\frac{\left(1-\left(1-\frac{y}{2n}\right)^{2n}(y+1)\right)^2}{\left(1-\frac{y}{2n}\right)^{2n}y^2}  \mathbb{1}_{0\leq y\leq x\leq 2n}.
\end{equation*}
To use the dominated convergence theorem, we want to upper-bound $F_n$ independently on $n$.
Note that, for all $ 0 \leq y\leq 2n $, \begin{align}
(1-y)\leq \left(1-\frac{y}{2n}\right)^{2n} \label{R2}\\
0\leq 1-\left(1-\frac{y}{2n}\right)^{2n}\left(y+1\right)\leq1\label{R1}\\
y\in[0,2n ]\mapsto  \left(1-\frac{y}{2n}\right)^{2n}  e^{y} \quad \text{ is decreasing} .\label{R3}
\end{align}
In particular
\begin{equation}\label{R4}
    \frac{\left(1-\left(1-\frac{y}{2n}\right)^{2n}(y+1)\right)^2}{\left(1-\frac{y}{2n}\right)^{2n}y^2}\leq \frac{\left(1-\left(1-{y}\right)(y+1)\right)^2}{\left(1-\frac{y}{2n}\right)^{2n}y^2}= \frac{y^2}{\left(1-\frac{y}{2n}\right)^{2n}}.
\end{equation}
Moreover $y\in[0,2n ] \mapsto \left(1-\frac{y}{2n}\right)^{2n}$ is decreasing, thus, for $y\leq x$,
    $F_n(x,y)\leq  \frac{(x+1)^2}{x^2}y^2$
which is integrable on $\{0\leq y\leq x\leq \epsilon\}$, for any $1>\epsilon>0$ set.
\\
From the increasing of $y\in[0,2n ]\mapsto\frac{y^2}{\left(1-\frac{y}{2n}\right)^{2n}}$, using \eqref{R4}, \eqref{R3} and \eqref{R2}, for $0\leq y\leq \epsilon \leq x$:
$$F_n(x,y)\leq \frac{\epsilon^2}{1-\epsilon}\frac{e^{-x}( x+1)^2}{x^2}$$
which is integrable on $\{0\leq y\leq \epsilon \leq x\}$.\\
We now take $\epsilon\leq y\leq x\leq 2n$. From \eqref{R1} and \eqref{R3}:
\begin{align*}
    F_n(x,y)&\leq \frac{\left(1-\frac{x}{2n}\right)^{2n}}{\left(1-\frac{y}{2n}\right)^{2n}}\frac{(x+1)^2}{x^2}\frac{1}{y^2} \leq\frac{e^{-x}(x+1)^2}{x^2}\frac{e^y}{y^2} \end{align*}
which is integrable on $\{\epsilon \leq y\leq x\}$ using Fubini.

We now turn to the a.s. limit $F$ of $F_n$:
\begin{equation*}
    F(x,y)=e^{y-x}\frac{\left(x+1\right)^2}{x^2}\frac{\left(1-e^{-y}(y+1)\right)^2}{y^2} \mathbb{1}_{0\leq y\leq x < \infty}\in L^1(\mathbb{R}^2) .
\end{equation*}
By the dominated convergence theorem, we then obtain
\[\Var(\tau_n)\underset{n\to+\infty}{\sim}\frac{2}{(2n+1)^2}\int F(x,y)dxdy.\]

\section{Cutoff for the sphere}\label{Sec3}
\begin{pro}\label{PR3}
On the sphere, the first covering time $\tau_n$ is satisfying:
\begin{align*}
    \mathbb{E}[\tau_n]& =\frac{1}{n-1}\sum\limits_{k=1}^{n-1}\frac{1}{k}=\frac{\ln(n-1)}{n-1}+\frac{\gamma}{n-1}+\frac{1}{2(n-1)^2}+o\left(\frac{1}{n^2}\right)\\
    \Var(\tau_n)&\underset{n\to +\infty}{\sim}\frac{1}{(n-1)^2}\frac{\pi^2}{6}
\end{align*}
where $\gamma$ is the Euler–Mascheroni constant.
\end{pro}
The proof is postponed to the two following subsections (Subsections \ref{Subsec5} and \ref{Subsec6}). As previously, we get the following:
\begin{theo}[For the sphere]\label{T3} Let $X^n$ be the (two-times accelerated) Brownian Motion in $ \mathbb{S}^n$, and $\gamma$ the Euler–Mascheroni constant.
Let $ c \in \mathbb{R} $ and $t_n = \frac{\ln(n-1)}{n-1} + \frac{c}{n-1}$: 
\begin{itemize}
\item if  $ c > \gamma$ then
  $\limsup\limits_{n \to \infty}  \fs(\cL(X^n_{t_n}),\mathcal{U}_{n}) \le \frac{\pi^2}{6 (c- \gamma)^2}   $
\item  if  $ c < \gamma$ then
  $\liminf\limits_{n \to \infty}  \fs(\cL(X^n_{t_n}),\mathcal{U}_{n}) \ge 1- \frac{\pi^2}{6 (c- \gamma)^2}   $.
\end{itemize}
Hence  $X^n$ has a cutoff in separation at time $\frac{\ln(n-1)}{n-1}  $ with   window $ \frac{1}{n-1}$, and we have the above control.
\end{theo}

\subsection{Computation of the mean, first part of the proof of Proposition \ref{PR3}}\label{Subsec5}
For the $n$-dimensional sphere, 
we remind from \eqref{EIn}:
\begin{equation}
    I_n(r)=\int_0^r \sin^{n-1}(x)dx \text{ and } I_n(\pi)=2W_{n-1}
\end{equation}
where $W_{n-1}$ denotes the corresponding Wallis' integral. As $I_n(\pi)-I_n(r)=I_n(\pi-r)$ for all $r\in[0,\pi]$:
\begin{equation}\label{ESphere}
    \EE[\tau_n]=\frac{1}{I_n(\pi)}\int_0^{\pi}\frac{I_n(x)I_n(\pi-x)}{\sin^{n-1}(x)}dx.
\end{equation}
Noticing that, for $n\geq 3$ \begin{align*}
I_n(r)&=\int_0^r \sin(x)\sin^{n-2}(x) \, dx\\&=[-\sin^{n-2}(x)\cos(x)]_0^r+(n-2)\int_0^r\sin^{n-3}(x)\cos^2(x)dx\\
&=-\sin^{n-2}(r)\cos(r)+(n-2)I_{n-2}(r)-(n-2)I_n(r),
\end{align*} we obtain the following recursive relations for $(I_n)_n$:
\begin{lem}\label{L2}
For all $n\geq 3$:
\begin{align*}
    I_n(r)&=\frac{1}{n-1}\left(-\sin^{n-2}(r)\cos(r)+(n-2)I_{n-2}(r)\right)\\
    I_n(\pi-r)&=\frac{1}{n-1}\left(\sin^{n-2}(r)\cos(r)+(n-2)I_{n-2}(\pi-r)\right).
\end{align*}In particular, we retrieve the recursive relation for Wallis' integrals: $$I_n(\pi)=\frac{n-2}{n-1}I_{n-2}(\pi).$$
\end{lem}
By obtaining a recursive relation for $\left((n-1)\EE[\tau_n]\right)_n$, we obtain an expression for $\EE[\tau_n]$ in function of the harmonic sums which proves the first part of Proposition \ref{PR3}:
\begin{pro}\label{P3}
For $n\geq 3$, $$(n-1)\EE[\tau_n]=(n-3)\EE[\tau_{n-2}]+\frac{1}{n-1}+\frac{1}{n-2}.$$
In particular $\EE[\tau_2]=1$ and $\EE[\tau_3]=\frac{3}{4}$.

\end{pro}
\proof
For $n\geq 3$, we have:
\begin{align}
I_n(\pi)\EE[\tau_n]&=\int_0^\pi\frac{1}{\sin^2(x)}\frac{I_n(x)I_n(\pi-x)}{\sin^{n-3}(x)}dx\notag\\
&=\left[-\cot (x)\frac{I_n(x)I_n(\pi-x)}{\sin^{n-3}(x)}\right]_0^{\pi}\label{E1}\\
&+\int_0^{\pi}\cot(x)\frac{\sin^{n-1}(x)\left(I_n(\pi-x)-I_n(x)\right)}{\sin^{n-3}(x)}dx\label{E2}\\
&-(n-3)\int_0^{\pi}\cot(x)\frac{\cos(x)}{\sin^{n-2}(x)}I_n(x)I_n(\pi-x)dx\label{E3}.
\end{align}
As $\left|\frac{I_n(x)}{\sin^{n-2}(x)}\right|\leq x\sin(x)\xrightarrow[x\to 0]{}0$, by symmetry of the sinus function, \eqref{E1} vanishes.
For \eqref{E2}, we get:
\begin{align*}\int_0^{\pi}\cot(x)&\frac{\sin^{n-1}(x)\left(I_n(\pi-x)-I_n(x)\right)}{\sin^{n-3}(x)}dx\\
&=\int_0^{\pi}\cos(x)\sin(x)\left(I_n(\pi-x)-I_n(x)\right)dx\\
&=\frac{1}{2}\left[\sin^2(x)\left(I_n(\pi-x)-I_n(x)\right)\right]_0^{\pi}+\int_0^\pi \sin^{n+1}(x)dx\\
&=I_{n+2}(\pi)=\frac{n}{n+1}I_n(\pi).
\end{align*}
We now turn to \eqref{E3}. We have: \begin{align*}
\int_0^{\pi}\cot(x)\frac{\cos(x)}{\sin^{n-2}(x)}I_n(x)I_n(\pi-x)dx
&=I_n(\pi)\EE[\tau_n]-\int_0^{\pi}\frac{I_n(x)I_n(\pi-x)}{\sin^{n-3}(x)}dx.
\end{align*}
Applying Lemma \ref{L2} on $I_n(x)I_n(\pi-x)$:
\begin{align*}
\int_0^{\pi}&\frac{I_n(x)I_n(\pi-x)}{\sin^{n-3}(x)}dx\\
&=\left(\frac{n-2}{n-1}\right)^2\int_0^{\pi}\frac{I_{n-2}(x)I_{n-2}(\pi-x)}{\sin^{n-3}(x)}dx-\frac{1}{(n-1)^2}\int_0^{\pi}\frac{\cos^2(x)\sin^{2(n-2)}(x)}{\sin^{n-3}(x)}dx\\
&+\frac{n-2}{(n-1)^2}\int_0^{\pi}\frac{\sin^{n-2}(x)\cos(x)\left(I_{n-2}(x)-I_{n-2}(\pi-x)\right)}{\sin^{n-3}(x)}dx\\
&=\left(\frac{n-2}{n-1}\right)^2I_{n-2}(\pi)\EE[\tau_{n-2}]-\frac{1}{(n-1)^2}\int_0^{\pi}(1-\sin^2(x))\sin^{n-1}(x)dx\\
&+\frac{n-2}{(n-1)^2}\int_0^{\pi}\sin(x)\cos(x)\left(I_{n-2}(x)-I_{n-2}(\pi-x)\right)dx\\
&=\frac{n-2}{n-1}I_{n}(\pi)\EE[\tau_{n-2}]-\frac{1}{(n-1)^2}\left(I_{n}(\pi)-I_{n+2}(\pi)\right)\\
&+\frac{n-2}{(n-1)^2}\left(\left[\frac{\sin^2(x)}{2}\left(I_{n-2}(x)-I_{n-2}(\pi-x)\right)\right]_0^{\pi}-\int_0^{\pi}\sin^{n-1}(x)dx\right)\\
&=I_n(\pi)\left(\frac{n-2}{n-1}\EE[\tau_{n-2}]-\frac{1}{(n-1)^2}\left(1-\frac{n}{n+1}\right)-\frac{n-2}{(n-1)^2}\right).
\end{align*}
Thus \eqref{E3} is equal to $(n-3)I_n(\pi)\left(-\EE[\tau_n]+\frac{n-2}{n-1}\EE[\tau_{n-2}]-\frac{1}{(n-1)^2(n+1)}-\frac{n-2}{(n-1)^2}\right)$.
Combining \eqref{E1},\eqref{E2} and \eqref{E3}, we obtain:
\begin{align*}
(n-2)\EE[\tau_n]&=\frac{n}{n+1}+\frac{(n-2)(n-3)}{n-1}\EE[\tau_{n-2}]-\frac{n-3}{(n-1)^2(n+1)}-\frac{(n-2)(n-3)}{(n-1)^2}\\
&=\frac{(n-2)(n-3)}{n-1}\EE[\tau_{n-2}]+\frac{2n-3}{(n-1)^2}\\
&=\frac{n-2}{n-1}\left((n-3)\EE[\tau_{n-2}]+\frac{1}{n-1}+\frac{1}{n-2}\right).
\end{align*}
This proves the recursive relation.
We now compute the first terms. As
   $ I_2(r)=\int_0^r \sin(x)dx=1-\cos(r)$,
\begin{align*}
    \EE[\tau_2]&=\frac{1}{I_2(\pi)}\int_0^{\pi} \frac{(1-\cos(r))(1+\cos(r))}{\sin(r)}dr=\frac{1}{I_2(\pi)}\int_0^{\pi} \frac{1-\cos^2(r)}{\sin(r)}dr\\
    &=\frac{1}{I_2(\pi)}\int_0^{\pi} \sin(r)dr=1
    \end{align*}
    We also obtain
    \begin{align*}
    \EE[\tau_3]&=\frac{1}{I_3(\pi)}\int_0^{\pi}\frac{I_3(x)I_3(\pi-x)}{\sin^2(x)}dx\\
    &=\frac{1}{I_3(\pi)}\left(\left[-\cot(x)I_3(x)I_3(\pi-x)\right]_0^{\pi}+\int_0^{\pi}\cot(x)\sin^2(x)(I_3(\pi-x)-I_3(x))dx\right)\\
    &=\frac{1}{I_3(\pi)}\left(\left[\frac{\sin^2(x)}{2}(I_3(\pi-x)-I_3(x))\right]_0^{\pi}+\int_0^{\pi}\sin^4(x)dx\right)\\
    &=\frac{I_5(\pi)}{I_3(\pi)}=\frac{3}{4}.
\end{align*}
\wwtbp

\subsection{The computation of the variance, second part of the proof of Proposition \ref{PR3}}\label{Subsec6}
Using \eqref{E0} and the symmetry of the volume form density we have:
\begin{align*}
   \Var(\tau_n)&=\frac{2}{I_n(\pi)^2} \int_0^{\pi}\frac{{I_n}(u)^2}{{I_n}'(u)} \int_u^{\pi}\frac{I_n(\pi-s)^2}{{I_n}'(s)}  ds \, du\\
   &=\frac{2}{I_n(\pi)^2}\left( \int_0^{\frac{\pi}{2}}\frac{{I_n}(u)^2}{{I_n}'(u)} \left(\int_u^{\frac{\pi}{2}}\frac{I_n(\pi-s)^2}{{I_n}'(s)}ds+\int_{\frac{\pi}{2}}^{\pi}\frac{I_n(\pi-s)^2}{{I_n}'(s)}ds\right) du\right.\\
   &+\left.\int_{\frac{\pi}{2}}^{\pi}\frac{{I_n}(u)^2}{{I_n}'(u)} \int_{u}^{\pi}\frac{I_n(\pi-s)^2}{{I_n}'(s)}  ds \, du\right)\\
      &=\frac{2}{I_n(\pi)^2}\left( \int_0^{\frac{\pi}{2}}\frac{{I_n}(u)^2}{{I_n}'(u)} \left(\int_u^{\frac{\pi}{2}}\frac{I_n(\pi-s)^2}{{I_n}'(s)}ds+\int_0^{\frac{\pi}{2}}\frac{I_n(s)^2}{{I_n}'(s)}ds\right) du\right.\\
   &+\left.\int_0^{\frac{\pi}{2}}\frac{{I_n}(\pi-u)^2}{{I_n}'(u)} \int_0^u\frac{I_n(s)^2}{{I_n}'(s)}  ds \, du\right)\\
         &=\frac{2}{{I_n(\pi)}^2}\left( 2\int_0^{\frac{\pi}{2}}\frac{{I_n(u)}^2}{{I_n}'(u)} \int_u^{\frac{\pi}{2}}\frac{I_n(\pi-s)^2}{{I_n}'(s)}ds\, du+\left(\int_0^{\frac{\pi}{2}}\frac{I_n(u)^2}{{I_n}'(u)}du\right)^2\right).
\end{align*}
Developing $I_n(\pi-s)^2=\left(I_n(\pi)-I_n(s)\right)^2$, we get:
\begin{align*}
    \Var(\tau_n)&=\frac{4}{I_n(\pi)^2}\int_0^{\frac{\pi}{2}}\frac{{I_n}(u)^2}{{I_n}'(u)} \int_u^{\frac{\pi}{2}}\frac{I_n(\pi)\left(I_n(\pi)-2I_n(s)\right)}{{I_n}'(s)}ds\, du\\
    &+4\left(\frac{1}{I_n(\pi)}\int_0^{\frac{\pi}{2}}\frac{{I_n(u)}^2}{{I_n}'(u)}du\right)^2\\
    &=8K_n+4\tilde K_n^2
\end{align*}
where
\begin{equation}
    K_n:=\frac{1}{I_n(\pi)}\int_0^{\frac{\pi}{2}}\frac{I_n(u)^2}{{I_n}'(u)} \int_u^{\frac{\pi}{2}}\frac{I_n\left(\frac{\pi}{2}\right)-I_n(s)}{{I_n}'(s)}ds\, du
\end{equation}
and 
\begin{equation}
    \tilde K_n:=\frac{1}{I_n(\pi)}\int_0^{\frac{\pi}{2}}\frac{{I_n(u)}^2}{{I_n}'(u)} du.
\end{equation}
We now look for the asymptotic behavior of $K_n$ and $\tilde K_n$, using a reparametrization and an asymptotic expansion to make appear the equivalents $\frac{1}{n^2}$ and $\frac{1}{n}$ respectively.
We first look for $K_n$:
\begin{equation*}
    K_{n}=\frac{1}{I_n(\pi)}\int\limits_{\substack{  \\0\leq u_1,u_2\leq u\leq s\leq r\leq \frac{\pi}{2}}}\frac{ \sin^{n-1}(u_1)\sin^{n-1}(u_2) \sin^{n-1}(r)}{\sin^{n-1}(u)\sin^{n-1}(s)}du\, ds\,du_1\,du_2\,dr.
\end{equation*}
Use now the following change of variables
\begin{align*}
    \sin(u_i)&=\left(1-\frac{h_i}{n-1}\right)\sin(u),\,\, i\in\{1,2\}\\
    \sin(u)&=\left(1-\frac{h_3}{n-1} \right)\sin(s)\\
    \sin(s)&=\left(1-\frac{h_4}{n-1}\right)\sin(r).
\end{align*}
In particular, as the Jacobian matrix is triangular, the corresponding Jacobian determinant is $J_n(h_1,h_2,h_3,h_4,r)=\alpha_1\alpha_2\alpha_3\alpha_4$ with
\begin{align*}
    \alpha_i&=-\frac{\left(1-\frac{h_3}{n-1} \right)\left(1-\frac{h_4}{n-1} \right)\sin(r)}{(n-1)\sqrt{\cos^2(r)+\frac{\delta_n(h_i,h_3,h_4)}{n-1}\sin^2(r)}}\\
    &=-\frac{\left(1-\frac{h_3}{n-1} \right)\left(1-\frac{h_4}{n-1} \right)}{\sqrt{n-1}\sqrt{(n-1)\cot^2(r)+\delta_n(h_i,h_3,h_4)}},\,\, i\in\{1,2\}\\
    &\\
     \alpha_3
    &=-\frac{\left(1-\frac{h_4}{n-1} \right)}{\sqrt{n-1}\sqrt{(n-1)\cot^2(r)+\delta_n(0,h_3,h_4)}},\\
   \alpha_4
    &=-\frac{1}{\sqrt{n-1}\sqrt{(n-1)\cot^2(r)+\delta_n(0,0,h_4)}}
\end{align*}
where $\delta_n(a_1,a_2,a_3):=(n-1)\left(1-\left(1-\frac{a_1}{n-1}\right)^2\left(1-\frac{a_2}{n-1}\right)^2\left(1-\frac{a_3}{n-1}\right)^2\right)$. 
We get:
\begin{align*}
    K_{n}=\frac{1}{I_n(\pi)}\int\limits_{\substack{ 0\leq h_1,h_2,h_3,h_4\leq n-1\\
    0\leq r\leq \frac{\pi}{2}}}& \left(1-\frac{h_1}{n-1}\right) ^{n-1}\left(1-\frac{h_2}{n-1}\right) ^{n-1}\left(1-\frac{h_3}{n-1}\right) ^{n-1}\\
    &\times \sin^{n-1}(r)J_n(h_1,h_2,h_3,h_4,r)dh_1\,dh_2\,dh_3\,dh_4\,dr.
\end{align*}
We also make the change of variable $x=\sqrt{n-1}\cot(r)$. In particular
we have
$\sin(r)=\left(\frac{x^2}{n-1}+1\right)^{-\frac{1}{2}}$ and $dx=-\sqrt{n-1}\left(\frac{x^2}{n-1}+1\right)dr$. Thus:
\begin{equation*}
    K_{n}=\frac{1}{(n-1)^{2+\frac{1}{2}}}\frac{1}{I_n(\pi)}\int\limits_{\substack{ 0\leq h_1,h_2,h_3,h_4\leq n-1\\
    0\leq x\leq +\infty}}F_n(h_1,h_2,h_3,h_4,x)dh_1\,dh_2\,dh_3\,dh_4\,dx
\end{equation*}
where
\begin{align*}F_n(h_1,h_2,h_3,h_4,x)=&\frac{ \left(1-\frac{h_1}{n-1}\right) ^{n-1}\left(1-\frac{h_2}{n-1}\right) ^{n-1}\left(1-\frac{h_3}{n-1}\right) ^{n+1}}{\sqrt{x^2+\delta_n(h_1,h_3,h_4)}\sqrt{x^2+\delta_n(h_2,h_3,h_4)}\sqrt{x^2+\delta_n(0,h_3,h_4)}}\\
&\frac{\left(1-\frac{h_4}{n-1}\right) ^{3}\left(1+\frac{x^2}{n-1}\right)^{-\frac{(n-1)}{2}-1}}{\sqrt{x^2+\delta_n(0,0,h_4)}}.\end{align*}

Using equivalences for the Wallis integral, we have $I_n(\pi)\underset{n\to +\infty}{\sim}\sqrt{\frac{2\pi}{{n-1}}}$.
Note that, for $a_1,a_2,a_3\in[0,n-1]$, $\delta_n(a_1,a_2,a_3)\geq \max(a_1,a_2,a_3)$.
For $n\geq 3$, we also have 
\begin{equation}\label{BornesupPuissance}
\left(1+\frac{x^2}{n-1}\right)^{n-1}\geq 1+x^2+\frac{n-2}{2(n-1)}x^4\geq \left(1+\frac{x^2}{2}\right)^2.
\end{equation}
Thus
$F_n(h_1,h_2,h_3,h_4,x)$ can be upper bounded by \[\frac{e^{-h_1}e^{-h_2}e^{-h_3}}{\sqrt{h_1}\sqrt{h_2}\sqrt{h_3}\sqrt{h_4}\left(1+\frac{x^2}{2}\right)}\] which is integrable for $h_4\leq \epsilon$ with $\epsilon>0$ set. It can also be upper bounded by
\[\frac{e^{-h_1}e^{-h_2}e^{-h_3}}{h_4^2\left(1+\frac{x^2}{2}\right)}\]
which is integrable for $h_4\geq \epsilon$.
Moreover, $\delta_n(a_1,a_2,a_3)\xrightarrow[n\to +\infty]{}2(a_1+a_2+a_3)$.
Using dominated convergence, we obtain:
\begin{align*}
K_{n}&\underset{n\to +\infty}{\sim}\frac{1}{(n-1)^{2}\sqrt{{2}{\pi}}}\int\limits_{[0,+\infty[^{5}}\frac{{e^{-h_1}e^{-h_2}e^{-h_3}e^{-\frac{x^2}{2}}}}{\sqrt{x^2+2(h_1+h_3+h_4)}\sqrt{x^2+2(h_2+h_3+h_4)}}\\
&\qquad\qquad\qquad\qquad\qquad\qquad\times\frac{dh_1\,dh_2\,dh_3\,dh_4\,dx}{\sqrt{x^2+2(h_3+h_4)}\sqrt{x^2+2h_4}}\\
&\underset{n\to +\infty}{\sim}\frac{K}{(n-1)^{2}}\end{align*}
with 
\[K:=\frac{1}{\sqrt{{2}{\pi}}}\int\limits_{\substack{0\leq x\leq \tilde h_4\leq \tilde h_3\leq \tilde h_1, \tilde h_2<+\infty}
}e^{-\frac{\tilde h_1^2}{2}}e^{-\frac{\tilde h_2^2}{2}}e^{\frac{\tilde h_3^2}{2}}e^{\frac{\tilde h_4^2}{2}}e^{-\frac{x^2}{2}}d \tilde h_1\,d\tilde h_2\,d\tilde h_3\,d\tilde h_4\,dx\]
where we used the change of variables
\begin{align*}
 \tilde h_i^2&=x^2+2(h_i+h_3+h_4),\,\, i=1,2\\
    \tilde h_3^2&=x^2+2(h_3+h_4)\\
    \tilde h_4^2&=x^2+2h_4.
\end{align*}
For $\tilde{K}_n=\frac{1}{I_n(\pi)}\int\limits_{0\leq u_1,u_2\leq u\leq\frac{\pi}{2}}\frac{\sin^{n-1}(u_1)\sin^{n-1}(u_2)}{\sin(u)}du_1du_2du$, we use similar change of variables:
\begin{align*}
    \sin(u_i)&=\left(1-\frac{h_i}{n-1}\right)\sin(u),\,\, i\in\{1,2\}\\
    x&=\sqrt{n-1}\cot(u).
\end{align*}
We get:
\[\tilde K_n=\frac{1}{(n-1)^{1+\frac{1}{2}}}\frac{1}{I_n(\pi)}\int\limits_{\substack{0\leq h_1,h_2\leq n-1\\ 0\leq x\leq +\infty}}\tilde{F}_n(h_1,h_2,x)dh_1\,dh_2\,dx\]
with 
\begin{align*}\tilde F_n(h_1,h_2,x)&=\frac{ \left(1-\frac{h_1}{n-1}\right)^{n-1}\left(1-\frac{h_2}{n-1}\right) ^{n-1}\left(1+\frac{x^2}{n-1}\right)^{-\frac{(n-1)}{2}-1}}{\sqrt{x^2+\delta_n(h_1,0,0)}\sqrt{x^2+\delta_n(h_2,0,0)}}\\
&\leq\frac{e^{-h_1}e^{-h_2}}{\sqrt{h_1}\sqrt{h_2}\left(1+\frac{x^2}{2}\right)}\end{align*}
which is integrable. From the dominated convergence theorem, we obtain
\[\tilde K_n\underset{n\to+\infty}{\sim}\frac{1}{(n-1)\sqrt{{2}{\pi}}}\int\limits_{[0,+\infty[^3}\frac{e^{-h_1}e^{-h_2}e^{-\frac{x^2}{2}}}{\sqrt{x^2+2h_1}\sqrt{x^2+2h_2}}dh_1\,dh_2\,dx=\frac{\tilde K}{(n-1)}\]
with 
\[\tilde{K}=\frac{1}{\sqrt{{2}{\pi}}}\int\limits_{0\leq x\leq \tilde h_1,\tilde h_2< +\infty}e^{-\frac{\tilde h_1^2}{2}}e^{-\frac{\tilde h_2^2}{2}}e^{\frac{x^2}{2}}d\tilde h_1\,d\tilde h_2\,dx.\]
Finally, we obtain 
$$\Var(\tau_n) \underset{n\to+\infty}{\sim} \frac{1}{(n-1)^2}\left(8K+4\tilde K^2\right).$$
We now make the computations for $K$ and $\tilde K$. First notice that, using polar coordinates, for $x\geq 0$,
\begin{equation}
    \int\limits_{x\leq \tilde h_1,\tilde h_2}e^{-\frac{\tilde h_1^2+\tilde h_2^2}{2}}d\tilde h_1\,d\tilde h_2
    =2\int_0^{\frac{\pi}{4}}\int_{\frac{x}{\sin(\theta)}}^{\infty}re^{-\frac{r^2}{2}}drd\theta=2\int_0^{\frac{\pi}{4}}e^{-\frac{ x^2}{2\sin^2(\theta)}}d\theta. 
\end{equation}
Then 
\begin{align*}\tilde K&=\sqrt{\frac{2}{\pi}}\int_0^{\frac{\pi}{4}}\int_0^{+\infty}e^{-\frac{x^2}{2}\cot^2(\theta)}dxd\theta=\sqrt{\frac{2}{\pi}}\int_0^{\frac{\pi}{4}}\tan(\theta)d\theta\int_0^{+\infty}e^{-\frac{x^2}{2}}dx\\
&=\left[-\ln(\cos(\theta))\right]_0^{\frac{\pi}{4}}=\frac{\ln(2)}{2}.
\end{align*}
We begin the same way for computing $K$:
\begin{align*}K&=\sqrt{\frac{2}{\pi}}\int_0^{\frac{\pi}{4}}\int\limits_{0\leq x\leq\tilde h_4\leq\tilde h_3}e^{-\frac{\tilde h_3^2}{2}\cot^2(\theta)}e^{\frac{\tilde h_4^2}{2}}e^{-\frac{x^2}{2}}dx\,d\tilde h_4\,d\tilde h_3\,d\theta.
\end{align*}
Using the change of variables $\tilde h_4=\rho\sin(\phi)$, $\tilde h_3=\rho\cos(\phi)$, we get:
\begin{align*}
K&=\sqrt{\frac{2}{\pi}}\int_0^{\frac{\pi}{4}}\int_0^{+\infty}e^{-\frac{x^2}{2}}\int_0^{\frac{\pi}{4}}\int_{\frac{x}{\sin(\phi)}}^{+\infty}\rho e^{-\frac{\rho^2}{2}\left(\cos^2(\phi)\cot^2(\theta)-\sin^2(\phi)\right)}d\rho\, d\phi\, dx\,d\theta\\
&=\sqrt{\frac{2}{\pi}}\int_0^{\frac{\pi}{4}}\int_0^{\frac{\pi}{4}}\int_0^{+\infty}e^{-\frac{x^2}{2}}\frac{e^{-\frac{x^2}{2\sin^2(\phi)}\left(\cos^2(\phi)\cot^2(\theta)-\sin^2(\phi)\right)}}{\left(\cos^2(\phi)\cot^2(\theta)-\sin^2(\phi)\right)} dx\, d\phi\,d\theta\\
&=\sqrt{\frac{2}{\pi}}\int_0^{\frac{\pi}{4}}\int_0^{\frac{\pi}{4}}\int_0^{+\infty}\frac{e^{-\frac{x^2}{2\tan^2(\phi)\tan^2(\theta)}}}{\left(\cot^2(\phi)\cot^2(\theta)-1\right)\sin^2(\phi)} dx\, d\phi\,d\theta\\
&=\int_0^{\frac{\pi}{4}}\int_0^{\frac{\pi}{4}}\frac{\tan(\phi)\tan(\theta)}{\left(\cot^2(\phi)\cot^2(\theta)-1\right)\sin^2(\phi)}  d\phi\,d\theta\\
&=\int_0^{\frac{\pi}{4}}\int_0^{\frac{\pi}{4}}\frac{\tan(\phi)\tan^3(\theta)}{\left(1-\tan^2(\phi)\tan^2(\theta)\right)\cos^2(\phi)}  d\phi\,d\theta .
\end{align*}
Setting $y=\tan^2(\phi)\tan^2(\theta)$, we get:
\begin{align*}
    K&=\frac{1}{2}\int_0^{\frac{\pi}{4}}\int_0^{\tan^2(\theta)}\frac{\tan(\theta)}{1-y}dy\,d\theta=\frac{1}{2}\int_0^{\frac{\pi}{4}}\tan(\theta)\left[\ln(1-y)\right]_0^{\tan^2(\theta)}d\theta\\
    &=-\frac{1}{2}\int_0^{\frac{\pi}{4}}\tan(\theta)\ln(1-\tan^2(\theta))d\theta=-\frac{1}{4}\int_0^1 \frac{\ln(1-u)}{1+u}du\\
    &=\frac{1}{4}\int_0^1\sum\limits_{j\geq 1}\frac{u^j}{j}\sum\limits_{k\geq 0}(-u)^kdu=-\frac{1}{4}\sum\limits_{j\geq 1}\sum\limits_{k\geq 0}\frac{(-1)^j(-1)^{k+j+1}}{j(k+j+1)}
    \\
    &=-\frac{1}{8}\left(\sum\limits_{j,k\geq 1}\frac{(-1)^j(-1)^k}{jk}-\sum\limits_{j\geq 1}\frac{1}{j^2}\right)\\
    &=\frac{1}{8}\left(\frac{\pi^2}{6}-\ln(2)^2\right).
\end{align*}
Finally, we get:
\begin{equation*}
    \Var(\tau_n)\underset{n\to +\infty}{\sim}\frac{1}{(n-1)^2}\frac{\pi^2}{6}.
\end{equation*}

\section{ Cutoff for the real projective space}\label{Sec4}
\begin{pro}\label{PR4}
On the real projective space, the first covering time $\tau_n$ is satisfying:
\begin{align*}
    \mathbb{E}[\tau_n]& =\frac{2\ln(n-1)}{n-1}+\frac{2(\gamma-\ln(2))}{n-1}+o\left(\frac{1}{n}\right)\\
    \Var(\tau_n)&\underset{n\to +\infty}{\sim}\frac{4N}{(n-1)^2}
\end{align*}
where $\gamma$ is the Euler–Mascheroni constant and \begin{equation*}N=\frac{2^4}{\pi}\int_0^{+\infty}e^{\frac{\tilde h_4^2}{2}}\left(\int_0^{\tilde h_4}e^{\frac{-x^2}{2}}dx\right)^2\int_{\tilde h_4}^{+\infty}e^{\frac{{\tilde h_3}^2}{2}}\left(\int_{\tilde h_3}^{+\infty}e^{\frac{-{\tilde h_1}^2}{2}}d{\tilde h_1}\right)^2d{\tilde h_3}d{\tilde h_4}.\end{equation*}
\end{pro}
To prove this proposition, we use the closeness between the volume on the spheres and the real projective space. The proof is detailed in the two following subsections (Subsections \ref{Subsec7} and \ref{Subsec8}).
Again, following the proof of Theorem \ref{T1}, we obtain:
\begin{theo}[For the real projective space]\label{T4} Let $X^n$ be the (two-times accelerated) Brownian Motion in $ \mathbb{RP}^n$, and $\gamma$ the Euler–Mascheroni constant.
Let $ c \in \mathbb{R} $ and $t_n = \frac{2\ln(n-1)}{n-1} + \frac{2c}{n-1}$:
\begin{itemize}
\item if  $ c > \gamma-\ln(2)$ then
  $\limsup\limits_{n \to \infty}  \fs(\cL(X^n_{t_n}),\mathcal{U}_{n}) \le \frac{N}{(c-\left(\gamma-\ln(2)\right))^2}   $
\item  if  $ c < \gamma -\ln(2) $ then
  $\liminf\limits_{n \to \infty}  \fs(\cL(X^n_{t_n}),\mathcal{U}_{n}) \ge 1- \frac{N}{(c- \left(\gamma-\ln(2)\right))^2}   $.
\end{itemize}
Hence  $X^n$ has a cutoff in separation at time $\frac{2\ln(n-1)}{n-1}  $ with   window $ \frac{2}{n-1}$, and we have the above control.
\end{theo}

\subsection{Computation of the mean, first part of the proof of Proposition \ref{PR4}}\label{Subsec7}
For the real projective space $\mathbb{RP}^n$ ($a=1$), the expression for $I_n(r)$ given by~\eqref{EIn} is quite close to the one for the sphere. Some of the results obtained in Section \ref{Sec3} will then be useful to study the first covering time on $\mathbb{RP}^n$. In this section we will use the notations $I_n^{\mathbb{S}}$ and $\tau_n^{\mathbb{S}}$ when the quantities relate to the the sphere and keep the notations $I_n$ and $\tau_n$ when they relate to the real projective space. We have:  \begin{equation}I_n(r)=\int_0^{r}\sin^{n-1}\left(\frac{x}{2}\right)dx=2\int_0^{\frac{r}{2}}\sin^{n-1}\left(x\right)dx=2I_n^{\mathbb{S}}\left(\frac{r}{2}\right).
\end{equation}
In particular, $I_n(\pi)=2I_n^{\mathbb{S}}\left(\frac{\pi}{2}\right)=2W_{n-1}=I_n^{\mathbb{S}}(\pi)$.
Thus, \begin{align*}
    \EE[\tau_n]&=\frac{1}{I_n(\pi)}\int_0^{\pi}\frac{I_n(x)\left(I_n(\pi)-I_n(x)\right)}{I_n'(x)}dx\\
    &=\frac{2}{W_{n-1}}\int_0^{\pi}\frac{I_n^{\mathbb{S}}\left(\frac{x}{2}\right)\left(I_n^{\mathbb{S}}\left(\frac{\pi}{2}\right)-I_n^{\mathbb{S}}\left(\frac{x}{2}\right)\right)}{{I_n^{\mathbb{S}}}'\left(\frac{x}{2}\right)}dx\\
    &=\frac{4}{W_{n-1}}\int_0^{\frac{\pi}{2}}\frac{I_n^{\mathbb{S}}\left(x\right)\left(\frac{1}{2}I_n^{\mathbb{S}}\left(\pi\right)-I_n^{\mathbb{S}}\left(x\right)\right)}{{I_n^{\mathbb{S}}}'\left(x\right)}dx\\
    &=\frac{4}{W_{n-1}}\int_0^{\frac{\pi}{2}}\frac{I_n^{\mathbb{S}}\left(x\right)\left(I_n^{\mathbb{S}}\left(\pi\right)-I_n^{\mathbb{S}}\left(x\right)\right)}{{I_n^{\mathbb{S}}}'\left(x\right)}dx-4\int_0^{\frac{\pi}{2}}\frac{I_n^{\mathbb{S}}\left(x\right)}{{I_n^{\mathbb{S}}}'\left(x\right)}dx
\end{align*}
From \eqref{ESphere}, noticing that $$\int_{\frac{\pi}{2}}^{\pi}\frac{I_n^{\mathbb{S}}(x)I_n^{\mathbb{S}}(\pi-x)}{{I_n^{\mathbb{S}}}'(x)}dx=\int_0^{\frac{\pi}{2}}\frac{I_n^{\mathbb{S}}(\pi-x)I_n^{\mathbb{S}}(x)}{{I_n^{\mathbb{S}}}'(x)}dx,$$ the mean for the sphere is given by:
\begin{equation*}
    \EE[\tau_n^{\mathbb{S}}]=\frac{1}{W_{n-1}}\int_0^{\frac{\pi}{2}}\frac{I_n^{\mathbb{S}}(x)\left(I_n^{\mathbb{S}}\left(\pi\right)-I_n^{\mathbb{S}}(x)\right)}{{I_n^{\mathbb{S}}}'(x)}dx.
\end{equation*}
Thus, for the projective space, we have:
 \begin{align*}
    \EE[\tau_n]&=4\EE[\tau_n^{\mathbb{S}}]-4E_n.
\end{align*}
where
$E_n:=\int_0^{\frac{\pi}{2}}\frac{I_n^{\mathbb{S}}\left(x\right)}{{I_n^{\mathbb{S}}}'\left(x\right)}dx$. We are now left to compute $E_n$. As for the computation of the mean for the sphere we use a recursive relation:
\begin{lem}\label{L4}
For $n\geq 3$, we have
$$(n-1)E_n=\frac{1}{n-2}+(n-3)E_{n-2}$$
and $$E_2=\ln(2)\,, E_3=\frac{1}{2}.$$
In particular, \begin{equation}\label{E_nReal}
E_n=\frac{1}{n-1}\sum\limits_{\substack{k=1\\
k \text{ odd}}}^{n-1}\frac{1}{k}+\frac{R_n}{n-1}\end{equation}
where 
\begin{equation*}R_n=\begin{cases}
\sum\limits_{k\geq n}\frac{(-1)^{k+1}}{k}=O\left(\frac{1}{n-1}\right)\text{ if }n\text{ even}\\
0\text{ if }n\text{ odd}
\end{cases}\end{equation*}
\end{lem}
From Lemma \ref{L4} and Proposition \ref{P3}, we then obtain:
\begin{lem}\label{P5}
For the real projective space of dimension $n$, the mean satisfies 
$$\EE[\tau_n]=
\frac{4}{n-1}\sum\limits_{k=1}^{\left\lfloor{\frac{n-1}{2}}\right\rfloor}\frac{1}{2k}-4\frac{R_n}{n-1}.$$

In particular, this proves the first part of Proposition \ref{PR4}.
\end{lem}
\prooff{Proof of Lemma \ref{L4}}
Set $n\geq 3$. 
Using integration by part as for Proposition \ref{P3}, we get:
\begin{align*}
    E_n&=\int_0^{\frac{\pi}{2}}\frac{1}{\sin^2(r)}\frac{I_n^{\mathbb{S}}(r)}{\sin^{n-3}(r)}dr\\
    &=\left[-\cot(r)\frac{I_n^{\mathbb{S}}(r)}{\sin^{n-3}(r)}\right]_0^{\frac{\pi}{2}}-(n-3)\int_0^{\frac{\pi}{2}}\cot(r)\frac{\cos(r)}{\sin^{n-2}(r)}I_n^{\mathbb{S}}(r)dr\\
    &+\int_0^{\frac{\pi}{2}}\frac{\cot(r)}{\sin^{n-3}(r)}\sin^{n-1}(r)dr.
\end{align*}
As for \eqref{E1}, the first term vanishes. Thus:
\begin{align*}
     E_n&=-(n-3)\int_0^{\frac{\pi}{2}}\frac{\cos^2(r)}{\sin^{n-1}(r)}I_n^{\mathbb{S}}(r)dr+\int_0^{\frac{\pi}{2}}\cos(r)\sin(r)dr\\
     &=-(n-3)\int_0^{\frac{\pi}{2}}\frac{1-\sin^2(r)}{\sin^{n-1}(r)}I_n^{\mathbb{S}}(r)dr+\left[\frac{\sin^2(r)}{2}\right]_0^{\frac{\pi}{2}}\\
    &=-(n-3)E_n +(n-3)\int_0^{\frac{\pi}{2}}\frac{I_n^{\mathbb{S}}(r)}{\sin^{n-3}(r)}dr+\frac{1}{2}.
\end{align*}
Gathering the terms with $E_n$ and using Lemma \ref{L2}, we get:
\begin{align*}
    (n-2)E_n&=-\frac{n-3}{n-1}\int_0^{\frac{\pi}{2}}\frac{\sin^{n-2}(r)\cos(r)}{\sin^{n-3}(r)}dr+\frac{(n-2)(n-3)}{n-1}E_{n-2}+\frac{1}{2}\\
    &=-\frac{n-3}{n-1}\left[\frac{\sin^2(r)}{2}\right]_0^{\frac{\pi}{2}}+\frac{(n-2)(n-3)}{n-1}E_{n-2}+\frac{1}{2}\\
    &=-\frac{n-3}{2(n-1)}+\frac{(n-2)(n-3)}{n-1}E_{n-2}+\frac{1}{2}\\
      &=\frac{1}{n-1}+\frac{(n-2)(n-3)}{n-1}E_{n-2}.
\end{align*}
This provides the expected recursive relation. In particular we get:
$$(n-1)E_n=\begin{cases}
\sum\limits_{\substack{2\leq k\leq n-1\\ k \text{ even}}}\frac{1}{k}+E_2\text{ if }n\text{ even}\\
\sum\limits_{\substack{2\leq k\leq n-1\\ k\text{ odd }}}\frac{1}{k}+2E_3\text{ if }n\text{ odd}
\end{cases}.$$
We now turn to the computation of the first terms:
\begin{align*}
    E_2&=\int_0^{\frac{\pi}{2}}\frac{\int_0^r\sin(s)ds}{\sin(r)}dr=\int_0^{\frac{\pi}{2}}\frac{1-\cos(r)}{\sin(r)}dr=\int_0^{\frac{\pi}{2}}\frac{\sin\left(\frac{r}{2}\right)}{\cos\left(\frac{r}{2}\right)}dr\\
    &=\left[-2\ln\left(\cos\left(\frac{r}{2}\right)\right)\right]_0^{\frac{\pi}{2}}=\ln(2).
    \end{align*}
    In particular, using the series expression $\ln(2)=\sum\limits_{k\geq 1}\frac{(-1)^{k+1}}{k}$ and noticing that $$\left|\sum\limits_{k\geq n}\frac{(-1)^{k+1}}{k}\right|\leq \frac{1}{n},$$we obtain \eqref{E_nReal} for $n$ even.
For $E_3$, still using the fact that $\frac{I_n^{\mathbb{S}}(r)}{\sin^{n-2}(r)}\xrightarrow[r\to 0]{}0$ we get:
\begin{align*}
    E_3&=\int_0^{\frac{\pi}{2}}\frac{\int_0^r\sin^2(s)ds}{\sin^2(r)}dr=\left[-\cot(r)\int_0^r\sin^2(s)ds\right]_0^{\frac{\pi}{2}}+\int_0^{\frac{\pi}{2}}\cot(r)\sin^2(r)dr\\
    &=\int_0^{\frac{\pi}{2}}\cos(r)\sin(r)dr=\frac{1}{2}.
\end{align*}
This gives \eqref{E_nReal} for $n$ odd.
\wwtbp

\subsection{The computation of the variance, second part of the proof of Proposition \ref{PR4}}\label{Subsec8}
From \eqref{E0}, we have:
\begin{align*}
\Var(\tau_n)&=\frac{8}{I_n(\pi)^2} \int_0^{\frac{\pi}{2}}\frac{I_n^2(2u)}{I'_n(2u)} \int_u^{\frac{\pi}{2}}\frac{(I_n(\pi)-I_n(2s))^2}{I'_n(2s)}  ds \, du\\
&=\frac{8}{I_n^{\mathbb{S}}(\pi)^2} \int_0^{\frac{\pi}{2}}\frac{4{I_n^{\mathbb{S}}}(u)^2}{{I^{\mathbb{S}}_n}'(u)} \int_u^{\frac{\pi}{2}}\frac{(2I_n^{\mathbb{S}}\left(\frac{\pi}{2}\right)-2I_n^{\mathbb{S}}(s))^2}{{I_n^{\mathbb{S}}}'(s)}  ds \, du\\
&=\frac{2^7}{I_n^{\mathbb{S}}(\pi)^2} \int_0^{\frac{\pi}{2}}\frac{{I_n^{\mathbb{S}}}(u)^2}{{I^{\mathbb{S}}_n}'(u)} \int_u^{\frac{\pi}{2}}\frac{(I_n^{\mathbb{S}}\left(\frac{\pi}{2}\right)-I_n^{\mathbb{S}}(s))^2}{{I_n^{\mathbb{S}}}'(s)}  ds \, du\\
&=\frac{2^8}{I_n^{\mathbb{S}}(\pi)^2}
\int\limits_{\substack{0\leq u_1,u_2\leq u\leq s\leq r_1\leq r_2\leq \frac{\pi}{2}}}\frac{ \sin^{n-1}(u_1)\sin^{n-1}(u_2) \sin^{n-1}(r_1)\sin^{n-1}(r_2)}{\sin^{n-1}(u)\sin^{n-1}(s)}\\
& \hspace{250pt} du\,ds\,du_1\,du_2\,dr_1\,dr_2.
\end{align*}
As for the computation of $K_n$ and $\tilde{K}_n$ in Subsection \ref{Subsec6}, we make the change of variables:
\begin{align*}
\sin(u_i)&=\left(1-\frac{h_i}{n-1}\right)\sin(u),\, i\in\{1,2\}\\
\sin(u)&=\left(1-\frac{h_3}{n-1}\right)\sin(s)\\
\sin(s)&=\left(1-\frac{h_4}{n-1}\right)\sin(r_1)\\
\sin(r_1)&=\left(1-\frac{h_5}{n-1}\right)\sin(r_2)\\
x&=\sqrt{n-1}\cot(r_2).
\end{align*}
Then, \begin{equation*}\Var(\tau_n)=\frac{2^8}{I_n^{\mathbb{S}}(\pi)^2(n-1)^3}\int\limits_{\substack{ 0\leq h_1,h_2,h_3,h_4,h_5\leq n-1\\
    0\leq x\leq +\infty}}G_n(h_1,h_2,h_3,h_4,h_5,x)dh_1\,dh_2\,dh_3\,dh_4\,dh_5\,dx\end{equation*}
    with 
    \begin{align*}G_n(h_1,h_2,h_3,h_4,h_5,x)&=\frac{ \left(1-\frac{h_1}{n-1}\right) ^{n-1}\left(1-\frac{h_2}{n-1}\right) ^{n-1}\left(1-\frac{h_3}{n-1}\right) ^{n+1}\left(1-\frac{h_4}{n-1}\right) ^{3}}{\sqrt{x^2+\delta_n(h_1,h_3,h_4,h_5)}\sqrt{x^2+\delta_n(h_2,h_3,h_4,h_5)}\sqrt{x^2+\delta_n(0,h_3,h_4,h_5)}}\\
    &\times\frac{\left(1-\frac{h_5}{n-1}\right) ^{n+3}\left(1+\frac{x^2}{n-1}\right)^{-(n-1)-1}}{\sqrt{x^2+\delta_n(0,0,h_4,h_5)}\sqrt{x^2+\delta_n(0,0,0,h_5)}}\end{align*}
    where this time \begin{equation*}\delta_n(a_1,a_2,a_3,a_4):=(n-1)\left(1-\left(1-\frac{a_1}{n-1}\right)^2\left(1-\frac{a_2}{n-1}\right)^2\left(1-\frac{a_3}{n-1}\right)^2\left(1-\frac{a_4}{n-1}\right)^2\right).\end{equation*}
    As previously, for any $\eps>0$, $n\geq 3$, $G_n$ can be upper bounded by an integrable function:
    \begin{equation*}
    G_n(h_1,h_2,h_3,h_4,h_5,x)\leq\frac{ e^{-h_1}e^{-h_2}e^{-h_3}e^{-h_5}}{\sqrt{h_1h_2h_3h_4h_5}\left(1+\frac{x^2}{2}\right)^2}\ind{0<h_4<\eps}+\frac{ e^{-h_1}e^{-h_2}e^{-h_3}e^{-h_5}}{h_4^2\sqrt{h_5}\left(1+\frac{x^2}{2}\right)^2}\ind{\eps<h_4}.
    \end{equation*}
    From dominated convergence and the previous equivalents for $I_n^{\mathbb{S}}(\pi)$, we obtain:
    \begin{align*}\Var(\tau_n)&\underset{n\to+\infty}{\sim}\frac{2^8}{2\pi (n-1)^2}\int\limits_{\substack{ 0\leq h_1,h_2,h_3,h_4,h_5,x\leq +\infty}}\frac{ e^{-h_1}e^{-h_2}e^{-h_3}e^{-h_5}e^{-x^2}}{\sqrt{x^2+2(h_1+h_3+h_4+h_5)}}\\
    &\times\frac{d h_1\, dh_2\, dh_3\, dh_4\, dh_5\,dx}{\sqrt{x^2+2(h_2+h_3+h_4+h_5)}\sqrt{x^2+2(h_3+h_4+h_5)}\sqrt{x^2+2(h_4+h_5)}\sqrt{x^2+2h_5}}\\
    &\underset{n\to+\infty}{\sim}\frac{N}{(n-1)^2}
    \end{align*}
    with \begin{align*}N&=\frac{2^8}{2\pi}\int\limits_{\substack{ 0\leq x\leq \tilde h_5\leq \tilde h_4\leq \tilde h_3\leq \tilde h_1,\tilde h_2\leq +\infty}} e^{-\frac{\tilde h_1^2}{2}}e^{-\frac{\tilde h_2^2}{2}}e^{\frac{\tilde h_3^2}{2}}e^{\frac{\tilde h_4^2}{2}}e^{-\frac{\tilde h_5^2}{2}}e^{-\frac{\tilde x^2}{2}}d\tilde h_1\,d\tilde h_2\,d\tilde h_3\,d\tilde h_4\,d\tilde h_5\,dx\\
    &=\frac{2^7}{2\pi}\int_0^{+\infty}e^{\frac{\tilde h_4^2}{2}}\left(\int_0^{\tilde h_4}e^{\frac{-x^2}{2}}dx\right)^2\int_{\tilde h_4}^{+\infty}e^{\frac{{\tilde h_3}^2}{2}}\left(\int_{\tilde h_3}^{+\infty}e^{\frac{-{\tilde h_1}^2}{2}}d{\tilde h_1}\right)^2d{\tilde h_3}d{\tilde h_4}.\end{align*}


\vskip2cm
\hskip70mm
\vbox{
\copy4
 \vskip5mm
 \copy5
  \vskip5mm
 \copy6
}

\end{document}